\documentclass[12pt]{article}
\usepackage{amsmath,amssymb,amsthm,booktabs}
\usepackage{mathrsfs,euscript}
\usepackage{geometry}
\usepackage{pifont}
\usepackage[OT1]{fontenc}
\usepackage[utf8]{inputenc}
\usepackage[colorlinks,citecolor=blue,urlcolor=blue]{hyperref}
\usepackage{txfonts}
\usepackage{indentfirst}
\usepackage{multirow}
\usepackage{float,subfig}

\usepackage{bm}
\usepackage{euscript}
\usepackage{graphicx}
\usepackage{multicol}
\usepackage[usenames,dvipsnames,svgnames,table]{xcolor}

\usepackage[numbers]{natbib}
\usepackage{a4wide}

\numberwithin{equation}{section} \theoremstyle{plain}
\newtheorem{theorem}{Theorem}[section]
\newtheorem{lemma}{Lemma}[section]
\newtheorem{corollary}{Corollary}[section]

\newtheorem{assumption}{Assumption}[section]

\newcommand\EN{\EuScript{N}}
\newcommand\lsigma{\underline{\sigma}^2}
\newcommand\usigma{\overline{\sigma}^2}
\newcommand\lmu{\underline{\mu}}
\newcommand\umu{\overline{\mu}}
\newcommand\gn{\EN_G(0,[\lsigma,\usigma])}

\makeatletter
\def\ps@pprintTitle{%
  \let\@oddhead\@empty
  \let\@evenhead\@empty
  \let\@oddfoot\@empty
  \let\@evenfoot\@oddfoot
}
\makeatother

\begin{document}

\newcommand\gai[1]{{\color{red}#1}}

\newcommand\tabfig[1]{\vskip5mm \centerline{\textsc{Insert #1 around here}}  \vskip5mm}

\vskip2cm

\title{Linear regression under model  uncertainty}

\author{Shuzhen Yang\thanks{Shandong University-Zhong Tai Securities
    Institute for Financial Studies, Shandong University, PR China,
    (yangsz@sdu.edu.cn). This work was supported by the National Key
    R\&D program of China (Grant No.2018YFA0703900), National Natural
    Science Foundation of China (Grant No.11701330), and Young
    Scholars Program of Shandong University.}\quad Jianfeng
  Yao\thanks{Department of Statistics and Actuarial Science, The
    University of Hong Kong,  Pokfulam Road,  Hong Kong SAR, Email:
    jeffyao@hku.hk}}

\date{\today}

\maketitle

\begin{abstract}
  We reexamine the classical linear regression model when the model is
  subject to two types of uncertainty: (i) some of covariates are either
  missing or completely inaccessible, and (ii) the variance of the
  measurement error is undetermined and changing according to a
  mechanism unknown to the statistician.  By following the recent theory
  of sublinear expectation, we propose to characterize  such mean and variance
  uncertainty in the response variable by two specific
  {\em nonlinear random variables}, which
  encompass an infinite family of
  probability distributions for the response variable in the sense of (linear) classical
  probability theory.
  The approach enables a family of estimators under various loss
  functions for the regression parameter and the parameters related to
  model uncertainty. The consistency of the estimators is established
  under mild conditions on the data generation process. Three
  applications are introduced to assess the quality of the approach
  including a forecasting model for the S\&P Index. 
  \end{abstract}

\bigskip


\section{Introduction}
\label{sec:intro}

 Robust regression has been actively developed  during the
years 1970-2000. A long catalogue of robust estimates for the regression coefficients
has appeared in the literature that includes  the $L_1$,
$M$,  $GM$, $RM$, $LMS$ and $LTS$, $S$, $MM$, $\tau$ and $SRC$
estimates among others.\footnote{Actually Huber complained that ``the collection of estimates
to choose from has become so extensive that it is worse than bewildering,
namely
counterproductive''. \citep[page 195]{Huber09}}
According to Huber, a robust procedure (or {\em stability}, see
\cite[page 5]{Huber09},
is ``in the sense that small deviations from the model
assumptions should impair the performance only slightly, that is, the latter
(described, say, in terms of the asymptotic variance of an estimate, or of the
level and power of a test) should be close to the nominal value calculated at
the model''.
The robust regression estimates above have been designed to
achieve such robustness while improving estimation efficiency and
protecting against unexpected  {\em procedure breakdown}.

Note that a central assumption  in this  robust statistics
literature is that the majority of the data under analysis
follows a distribution given by an assumed model.
Although the assumed  model can be very generic,
it however must be {\em unique} as requested by statistical theory in order to
enable inference about the model.
When the data under analysis significantly deviates from the assumed model,
inference runs out of the   set-up of traditional robust
statistics. Quoting again Huber,
``the interpretation of results obtained by blind robust estimators becomes
questionable when the fraction of contaminants is no longer small.''
\citep[page 198]{Huber09}

Originated from  the field of mathematical
finance, {\em model uncertainty} is a concept that can help
statisticians  deal with ``no longer small''  deviations of the data
from an assumed model in some precise contexts.
In an early work, \cite{GSCH89} proposed to tackle model ambiguity aversion by
the family of  max-min expected utility functions,
in a framework where data may follow
an  infinite family of models (or distributions).
The concept of model uncertainty and its applications in mathematical
finance are successively developed in the papers
\cite{ALP95,L95,CE02,C06,FS11}.  Particularly,
{\em coherent risk measures}  were
introduced in \cite{A99} to study both market risks and non-market risks.
Over the last  decade,
a fundamental concept of {\em sublinear expectation} was   developed in
\cite{Peng2004,Peng2005} which provides a general theory for
quantifying uncertainty about probability distribution of random
variables, and more generally, of stochastic processes.\footnote{In fact, the theory covers
  {nonlinear expectations} which are more general than the
  concept of sublinear expectations. However for  the purpose of this
  paper, it is sufficient to consider sublinear expectations.}
One important result of the theory is a  central
limit theorem (under sublinear expectation) that  bridges the general theory and
statistical data analysis under model or distribution
uncertainty. Parallel to the role of  a classical central limit theorem
to classical statistical inference,
a {\em nonlinear normal distribution} is introduced to approximate
asymptotic distributions of  large sums of independent variables.
This nonlinear normal distribution under
sublinear expectation is the celebrated {\em  G-normal distribution}.
This theory is fully developed in the recent  monograph
\cite{Peng2019}. (A review of relevant results in
Appendix~\ref{sec:prelim}).

This new theory of sublinear expectation leads to many questions to
explore in data analysis in situations where distribution
uncertainty is inherent to the data generation process under
consideration.  An example of such exploration is a recent work
\cite{Peng2020} where we constructed a new VaR predictor for financial
indexes which shows a significant advantage over  most of the
existing benchmark VaR predictors. A fundamental idea
underlying  \cite{Peng2020} is that, in parallel to classical
data analysis where the normal distribution is a natural choice  for
measurement errors or data fluctuations, the G-normal distribution can
serve as a primary tool for analyzing  data fluctuations when
distributions of data are subject to high uncertainty.  Such high
distribution uncertainty is indeed common in  financial
indexes such as the NASDAQ and S\&P 500 indexes.
The results obtained  in \cite{Peng2020} for VaR prediction
provide  a new confirmation of the existence of such distribution
uncertainty.
They  also  showcase the
power and usefulness of the new theory of sublinear expectation for
data analysis under model or distribution uncertainty.

In this paper we explore the implication of such model uncertainly in
the context of regression analysis.
Precisely, consider a 
$q$-dimensional deterministic covariate vector  $X\in \mathbb{R}^q$ and
a univariate dependent random variable $Y\in \mathbb{R}$ within a
regression model of the form
\begin{equation}\label{eq:line}
  Y=\beta^{\top} X+\eta+\varepsilon,
\end{equation}
where $\beta\in  \mathbb{R}^q$ is the vector of regression
coefficients. The novelty here is the terms $\eta$ and $\varepsilon$
which account for mean uncertainty, and variance uncertainty, respectively.
In layman's language, we can say that
$\beta^{\top} X$ accounts for the contribution to the response mean
from  the given covariates $X$,  while the unexplained or remaining part of the mean
is non-accessible either because no more significant covariates are
available, or it is varying through a somehow  unknown mechanism.
This uncertain part of the mean is
modeled by the {\em nonlinear random variable} $\eta$.
Furthermore, the fluctuation of $Y$ around its true mean, that is, the
error $\varepsilon$, cannot be
determined by a single classical probability distribution;
rather it will follow the {\em nonlinear G-normal distribution} in order to
capture the underlying  uncertainty.
The model~\eqref{eq:line} is referred as {\em distribution-uncertain
  regression model}.

Under both uncertainties  about the mean and variance of the response
variable,  is  it  still possible to
 consistently  estimate  the regression parameter $\beta$ in
 \eqref{eq:line}?
 To answer the question,   we consider a general loss function
$\phi$ and introduce two {\em population optimal parameters},  under model
uncertainty, namely, 
\begin{equation}
  \label{eq:optimal-beta}
  \overline{\beta}^*(\phi)=\arg \min_{\beta\in \mathbb{R}^q} \mathbb{E}[\phi(Y-\beta^{\top} X-\eta)],
\end{equation}
and
\begin{equation}
  \label{eq:optimal-beta1}
  \underline{\beta}^*(\phi)=\arg \min_{\beta\in \mathbb{R}^q} -\mathbb{E}[-\phi(Y-\beta^{\top} X-\eta)].
\end{equation}
The particular feature here is that $ \mathbb{E}$ is the sublinear
expectation operator.
Possible   choices for the loss  function are 
$\phi(z)=z^2$ for the square loss,  $\phi(z)=[\alpha-I(z<0)]z$ for quantile loss at a
given level $\alpha\in(0,1)$,
and $\phi(z)=I(z\le 0)$ for the Value-at-Risk (VaR) loss.
In general the optimal parameters
$\overline{\beta}^*(\phi)$ and $\underline{\beta}^*(\phi)$ depend  on the loss function $\phi(\cdot)$ under
model uncertainty. On the other hand,  if $Y$ had neither mean
uncertainty nor variance uncertainty,
that is, $\eta$ was a real constant and
 $\varepsilon$ a classical centred noise variable,
the model~\eqref{eq:line} would become a classical linear regression model,
and
we would  have
$\overline{\beta}^*(\phi)=\underline{\beta}^*(\phi)\equiv \beta$ for a
large class of possible loss functions $\phi$.

As a  main contribution of the paper,
Theorem~\ref{the:minmax} in Section~\ref{sec:rlr}
characterizes  the population optimal parameters $\overline{\beta}^*(\phi)$ and $\underline{\beta}^*(\phi)$
for a wide class of convex loss functions $\phi$.
Next, in Section~\ref{sec:lse}
we apply this characterization to the case of
the square loss $\phi(z)=z^2$. Based on  this characterization,
we propose a class of estimators for
both the regression parameter $\beta$ and those parameters that
involve  in the mean-uncertainty variable $\eta$ and
the variance-uncertainty variable $\varepsilon$.
Under appropriate conditions on the data observation process, we
establish  large sample consistency of these estimators.

The related  literature on regression analysis under model
uncertainty is actually quite limited.
When the error distribution in the regression model belongs to a
{\em finite family}, \cite{Lin16} constructed a $k$-sample maximum
expectation regression over the given finite family of distribution.
Using the square loss, several estimators
are proposed which are consistent  and asymptotically  normal. In a follow-up work, still under the
assumption of finite-number uncertainty, \cite{Lin17} investigated a
more general form of maximum expectation regression estimators and
established their consistency and asymptotic normality under
appropriate conditions.

Other sections of the paper are as follows.
Section~\ref{sec:simul} reports simulation  experiments to assess the
finite-sample properties of the proposed estimators under model 
uncertainty.
In Section~\ref{sec:app},
we develop three applications of our method to  robust regression,
regression under heteroscedastic error, and to an analysis of
daily returns of  the  S\&P 500 Index. 
In Appendix~\ref{sec:prelim}, we recall    useful results from the
theory of sublinear expectation which are
relevant to the work in this paper. All technical proofs are gathered
in Appendix~\ref{sec:proofs}.

%
\section{Linear regression under distribution uncertainty}
\label{sec:rlr}

Consider the distribution-uncertain regression model~\eqref{eq:line}
and the associate population optimal parameters $\overline{\beta}^*(\phi)$ and $\underline{\beta}^*(\phi)$ in
\eqref{eq:optimal-beta} and \eqref{eq:optimal-beta1}
for a given convex loss function $\phi$.
As mentioned in Introduction, standard choices for the loss function
cover the least squares estimator, quantile regression estimator and a VaR estimator.

Technically, we first construct a specific  infinite family of
probabilities. Consider a  canonical  probability space
$(\Omega,\mathcal{F},P)$ where
$\Omega=C([0,1])$ is the space of real-valued continuous functions on
[0,1].  Let  $\{B_t\}_{0\leq t\leq 1}$ be a Brownian motion. The parameter space we consider is
$\Theta=L^2(\Omega\times[0,1],[\underline{\sigma},\overline{\sigma}])$,
space of square-integrable, progressively measurable random processes
on [0,1] with values in
the interval $[\underline{\sigma},\overline{\sigma}]$.
Here the two parameters $0< \underline{\sigma}< \overline{\sigma}$ are
the lower and the upper limit for parameter processes
$\theta=(\theta_s)_{0\le t\le 1 }\in \Theta$, respectively.
The family of  probability measures $\{P_{\theta}\}$ is:   for $A\in\mathcal{F}$,
\[
P_{\theta}( A )=P \circ \xi_{\theta}^{-1}(A)=P(\xi_{\theta}\in A),
\quad
\text{where\quad}
\xi_{\theta}(\cdot)=\int_0^{\centerdot}\theta_sdB_s.
\]

The infinite family of probabilities
$\{P_{\theta}\}_{\theta\in\Theta}$  will govern the regression
model~\eqref{eq:line}. Precisely, under $P_{\theta}$,  the mean uncertainty
variable $\eta$ takes a constant $\mu_{\theta}\in[\lmu,\umu],\ \theta\in\Theta$,
while the variance uncertainty variable
$\varepsilon$ follows a nonlinear G-normal distribution $\gn$,
with lower and upper variance parameters
$(\lsigma,\usigma)$.\footnote{The details of G-normal distribution are given in Appendix \ref{sec:prelim-Gn}.}
Note that the distribution
uncertainty  of the error  $\varepsilon$ includes an infinite family
of distributions $\{F_{\theta}(\cdot)\}_{\theta\in\Theta}=\gn$, where
$F_{\theta}(\cdot)$ is determined  by $P_{\theta}$.


By the representation theorem of sublinear expectation,  Theorem~\ref{the:rep-sub}, we can express the nonlinear expectation of any
function of  $\varepsilon$ as
$$
\mathbb{E}[\phi(\varepsilon)]=\max_{\theta\in\Theta}E_{\theta}[\phi(\varepsilon)],
$$
where $E_{\theta}[\cdot]$ is the classical linear expectation under
$P_{\theta}$.
Therefore the population optimal parameters in
\eqref{eq:optimal-beta}-\eqref{eq:optimal-beta1}
have the form
\begin{equation}\label{eq:cos1}
\overline{\beta}^*(\phi)=\arg  \min_{\beta\in
  \mathbb{R}^q}\mathbb{E}[\phi(Y-\beta^{\top} X-\eta)]=\arg \min_{\beta\in \mathbb{R}^q}\max_{\theta\in \Theta}{E}_{\theta}[\phi(Y-\beta^{\top} X-\mu_{\theta})],
\end{equation}
and
\begin{equation}\label{eq:cos2}
\underline{\beta}^*(\phi)=\arg\min_{\beta\in
  \mathbb{R}^q}-\mathbb{E}[-\phi(Y-\beta^{\top} X-\eta)]=\arg \min_{\beta\in \mathbb{R}^q}\min_{\theta\in \Theta}{E}_{\theta}[\phi(Y-\beta^{\top} X-\mu_{\theta})].
\end{equation}
In other words,
 $\overline{\beta}^*(\phi)$ and $\underline{\beta}^*(\phi)$ are
optimal for the min-max loss and the min-min loss strategies,
respectively, over the infinite family of probabilities
$\{P_\theta\}_{\theta\in\Theta}$.

Next, we have a technical lemma of exchange rule between the
maximization or minimization steps in \eqref{eq:cos1} and \eqref{eq:cos2}.

\begin{lemma} \label{lem:minmax}
  We assume that the loss function $\phi(\cdot)\in C_{l.Lip}(\mathbb{R})$
is convex.  We have the exchange formulas for  (\ref{eq:cos1}) and (\ref{eq:cos2}):
  \begin{equation}\label{eq:minmax}
  \min_{\beta\in \mathbb{R}^q}\max_{\theta\in \Theta}{E}_{\theta}[\phi(Y-\beta^{\top} X-\mu_{\theta})]
    =\max_{\theta\in \Theta}\min_{\beta\in \mathbb{R}^q}{E}_{\theta}[\phi(Y-\beta^{\top} X-\mu_{\theta})],
  \end{equation}
  and
    \begin{equation}\label{eq:minmin}
  \min_{\beta\in \mathbb{R}^q}\min_{\theta\in \Theta}{E}_{\theta}[\phi(Y-\beta^{\top} X-\mu_{\theta})]
    =\min_{\theta\in \Theta}\min_{\beta\in \mathbb{R}^q}{E}_{\theta}[\phi(Y-\beta^{\top} X-\mu_{\theta})].
  \end{equation}
\end{lemma}

As a consequence of  Lemma~\ref{lem:minmax}, the  two  optimal parameters
$\overline{\beta}^*(\phi)$ and $\underline{\beta}^*(\phi)$ can
actually be determined under two classical normal distributions
$\EN(\mu_{\underline{\theta}_{\phi}},\lsigma)$ and
$\EN(\mu_{\overline{\theta}_{\phi}},\usigma)$, with some specific mean
parameters   $\mu_{\underline{\theta}_{\phi}  }$  and
$\mu_{\overline{\theta}_{\phi}}$.
This characterization of the parameters are instrumental for the
construction of their
estimators presented in Section~\ref{sec:lse}. 

\begin{theorem} \label{the:minmax}
We assume that the loss function $\phi(\cdot)\in C_{l.Lip}(\mathbb{R})
$ is convex. There
exists an optimal distribution parameter $\overline{\theta}_\phi(s)=\overline{\sigma},\ 0\leq s\leq 1$, such that
\[
\overline{\beta}^*(\phi)=\arg \min_{\beta\in \mathbb{R}^q}E_{\overline{\theta}_\phi}\left[ \phi(Y-\beta^{\top} X-\mu_{\overline{\theta}_\phi })\right].
\]
Similarly, there exists another optimal distribution parameter
$\underline{\theta}_\phi(s)=\underline{\sigma},\ 0\leq s\leq 1$, such that
\[
\underline{\beta}^*(\phi)=\arg \min_{\beta\in \mathbb{R}^q}E_{\underline{\theta}_\phi}\left[ \phi(Y-\beta^{\top} X-\mu_{\underline{\theta}_\phi })\right].
\]
\end{theorem}
The proofs of Lemma~\ref{lem:minmax} and
Theorem~\ref{the:minmax} are given in Appendix~\ref{proof-lem-minmax}
and \ref{proof-the-minmax}, respectively.

In  order to calculate  the two optimal parameters
$\overline{\beta}^*(\phi)$ and $\underline{\beta}^*(\phi)$,
we can use Theorem \ref{the:minmax} with the following  two-step procedure.

\begin{itemize}
\item[(1)\ ]
  Find the optimal liner expectations $E_{\overline{\theta}_{\phi}}[\cdot]$ and $E_{\underline{\theta}_{\phi}}[\cdot]$ based on the criterion function $\phi(\cdot)$ such that
  $$
  E_{\overline{\theta}_{\phi}}[\phi(Y-\beta^{\top} X-\mu_{\overline{\theta}_{\phi}})]=\max_{\theta\in \Theta}E_{\theta}[\phi(Y-\beta^{\top} X-\mu_{\theta})],
$$
and
$$
E_{\underline{\theta}_{\phi}}[\phi(Y-\beta^{\top} X-\mu_{\underline{\theta}_{\phi}})]=\min_{\theta\in \Theta}E_{\theta}[\phi(Y-\beta^{\top} X-\mu_{\theta})].
$$
\item[(2)\ ]
  Once $E_{\overline{\theta}_{\phi}}[\cdot]$ and
  $E_{\underline{\theta}_{\phi}}[\cdot]$ are found,  perform standard
  regression analysis under the two linear expectations to
  find the optimal parameters
$$
\overline{\beta}^*(\phi)=\arg \min_{\beta\in\mathbb{R}^q}E_{\overline{\theta}_{\phi}}[\phi(Y-\beta^{\top} X-\mu_{\overline{\theta}_{\phi}})],\quad \underline{\beta}^*(\phi)=\arg \min_{\beta\in\mathbb{R}^q}E_{\underline{\theta}_{\phi}}[\phi(Y-\beta^{\top} X-\mu_{\underline{\theta}_{\phi}})].
$$
\end{itemize} 
 
This two-step procedure defines a
new mechanism for  determining  the optimal parameters $\overline{\beta}^*(\phi)$ and $\underline{\beta}^*(\phi)$
under the considered distribution  uncertainty. The procedure is valid for a general convex loss  function $\phi(\cdot)\in C_{l.Lip}(\mathbb{R}^q)$.

\section{Least squares regression under distribution uncertainty}\label{sec:lse}

We now develop the least squares procedure for the estimation of the
regression parameter $\beta$ under the distribution-uncertain
model~(\ref{eq:line}).  The loss function is thus
$\phi(\cdot)=(\cdot)^2$, and the two population optimal parameters in
\eqref{eq:cos1} and   \eqref{eq:cos2} are:
\begin{equation}\label{eq:lse}
  \overline{\beta}^*=\arg \min_{\beta\in \mathbb{R}^q}\mathbb{E}[(Y-\beta^{\top}
    X-\eta)^2],\quad \underline{\beta}^*=\arg \min_{\beta\in \mathbb{R}^q}-\mathbb{E}[-(Y-\beta^{\top}
    X-\eta)^2].
\end{equation}
We call $\overline{\beta}^*$ the {\em upper-least squares parameter}
(U-LSE), and $\underline{\beta}^*$ the {\em lower-least squares parameter} (L-LSE).
Applying the general Theorem~\ref{the:minmax} to the present case, we
get the following characterization of these parameters, as well as
that of the two variance parameters $\underline\sigma$ and
$\overline\sigma$.

\begin{theorem}\label{the:u-lse}
Consider the distribution-uncertain regression model (\ref{eq:line})
under the square loss function $\phi(z)=z^2$.
\begin{enumerate}
\item[(i).\ ]
  The U-LSE $\overline{\beta}^*$ can be estimated by the observation samples from
  \begin{equation}\label{eq:u-lse}
    Y={\beta}^{\top}X+{\mu}_{\overline{\sigma}}+\varepsilon',
  \end{equation}
  where $\varepsilon'$ follows the  classical normal distribution $\EN(0,\overline{\sigma}^2)$.
\item[(ii).\ ]
  The L-LSE $\underline{\beta}^*$ can be estimated by the observation samples from
  \begin{equation}\label{eq:l-lse}
    Y={\beta}^{\top}X+{\mu}_{\underline{\sigma}}+\varepsilon'',
  \end{equation}
  where $\varepsilon''$ follows the  classical normal distribution $\EN(0,\underline{\sigma}^2)$.
\item[(iii).\ ]
  The variance parameters $\overline{\sigma}$ and $\underline{\sigma}$
  are characterized  as follows:
  \begin{equation}\label{eq:lse-sigma}
    \overline{\sigma}^2={E}_{\overline{\sigma}}[(Y-\overline{\beta}^* X-\mu_{\overline{\sigma}})^2],\quad \underline{\sigma}^2={E}_{\underline{\sigma}}[(Y-\underline{\beta}^* X-\mu_{\underline{\sigma}})^2],
  \end{equation}
  where $({E}_{\underline{\sigma}}[\cdot],{E}_{\overline{\sigma}}[\cdot])$ mean the expectations under $\overline{\theta}_{(\cdot)^2}(s)=\overline{\sigma},\ \underline{\theta}_{(\cdot)^2}(s)=\underline{\sigma},\ 0\leq s\leq 1$.
\end{enumerate}
\end{theorem}

Results in  Theorem \ref{the:u-lse} can be summarized as follows.
When we adopt  the min-max strategy,
$$
\min_{\beta\in \mathbb{R}^q}\max_{\theta\in\Theta}{E}_{\theta}[(Y-\beta X-\mu_{\theta})^2],
$$
the U-LSE $\overline{\beta}^*$ is the optimal parameter such that
\begin{equation}\label{eq:minmax0}
  \overline{\sigma}^2={E}_{\overline{\sigma}}[(Y-\overline{\beta}^* X-\mu_{\overline{\sigma}})^2]=\min_{\beta\in \mathbb{R}^q}\mathbb{E}[(Y-\beta^{\top}
X-\eta)^2].
\end{equation}
These characterizations will enable a sample counterpart of the U-LSE
$\overline{\beta}^*$
which will be a consistent estimator for the parameter $\beta$, and
subsequently, another  consistent estimator for the upper variance $\overline{\sigma}$.

Similarly, when we consider the min-min strategy, the L-LSE $\underline{\beta}^*$ is the optimal parameter such that
\begin{equation}\label{eq:minmin0}
  \underline{\sigma}^2={E}_{\underline{\sigma}}[(Y-\underline{\beta}^* X-\mu_{\underline{\sigma}})^2]=\min_{\beta\in \mathbb{R}^q}-\mathbb{E}[-(Y-\beta^{\top}
X-\eta)^2].
\end{equation}
Consistent estimators for both
the parameter $\beta$  and the lower variance  $\underline{\sigma}$
can also be derived by using the sample counterparts of these
parameters.

Consequently, we have the following results for U-LSE
$\overline{\beta}^*$ and L-LSE $\underline{\beta}^*$.
\begin{corollary}\label{coro:1}
For the given square loss function $\phi(z)=z^2$, we have that
$$
\overline{\beta}^*=\underline{\beta}^*=\beta
$$
for the distribution-uncertain regression model (\ref{eq:line}).
\end{corollary}

\subsection{Consistent estimators for the regression parameter $\beta$
  and distribution-uncertainty parameters  $(\underline\mu,\overline\mu,\underline\sigma^2,\overline\sigma^2)$}
\label{ssec:estim}

In order to formulate a theory of consistent estimation, we need to
define precisely the  generation process of the data under
consideration  as follows. 

\begin{description}
\item[Data generation process:] \quad
  The samples $\{(x_i,y_i)\}_{i=1}^T$ satisfy
  \begin{equation}\label{eq:DGP}
    y_i=\beta x_i+\eta_j+\varepsilon_i,\quad
    1+n_0(j-1)\leq i\leq n_0j,\ 1\leq j\leq K,
  \end{equation}
  where $\eta_j\in [\underline{\mu},\overline{\mu}]$, and $\varepsilon_i\sim  \EN(0,\sigma_j^2)$ with  $\sigma_j^2\in
  [\underline{\sigma}^2,\overline{\sigma}^2]$.
  Thus, there are $K$ groups in the samples, and each group has $n_0$
  elements with mean $\eta_j$ and variance $\sigma_j^2$. The total
  number of samples is $T=n_0K$.
\end{description}

The main challenge here for estimating the diverse parameters in the
distribution-uncertain model~(\ref{eq:line}) is that the theoretical characterizations of the U-LSE and L-LSE
parameters given in Theorem~\ref{the:u-lse}  cannot be used directly,
because we do not have at our disposal samples from the two normal
distributions
$\EN(\mu_{\underline{\sigma}},\lsigma)$
and
$\EN(\mu_{\overline{\sigma}},\usigma)$ that appear in \eqref{eq:u-lse}
and \eqref{eq:l-lse}, respectively. 
The difficulty is also due to the fact that from one sample
$(x_i,y_i)$ to next, the uncertain mean $\eta$ and uncertain error
$\varepsilon$ can change significantly. We propose a method based on
moving and overlapping blocks that lead to a family of {\em
  intermediate residuals} which are approximately distributed as
$\EN(\mu_{\underline{\sigma}},\lsigma)$.
These intermediate residuals are then used for consistent estimation
of  $(\beta,\lsigma$).
Afterwards, we can build consistent
estimations for $(\lmu,\umu,\usigma)$.

{Data generated under \eqref{eq:DGP} can be seen as a practical
  instance  of the general distribution-uncertain
  model~\eqref{eq:line}. It defines a specification needed for the
  introduction of an estimation theory. It is possible to relax a few
  conditions of the process. For example, 
  the group size $n_0$ may vary  with the groups, and the uncertain
  mean $\eta$ and uncertain error $\varepsilon$
  in the sample can have a controlled variation within  each group.
  Particularly,  only the  samples $\{(x_i,y_i)\}_{i=1}^T$ are available
to us, and we have no direct access to all  other parameters and
variables such as
(i) the group partition and the group length $n_0$;
(ii) the group means  $(\eta_j)$ that account for the mean uncertainty;
(iii) the error variances $(\sigma_j)$  that account for the error uncertainty.
Therefore, the problem of parameter estimation here is not
straightforward.

The main idea of our approach  is to use  {\em moving blocks}. The samples
$\{(x_i,y_i)\}_{i=1}^T$ are scanned subsequently as $m=T-n+1$ blocks
of a given block length $n$ as in
\[ \{1,\ldots,n\}, ~
   \{2,\ldots,n+1\},~ \ldots,~ \{T-n+1,\ldots,T\}.
\]
Denote the data in the $l$th block by  $B_l=\{(x_i,y_i)\}_{l\le i\le l+n-1}$,
$1\le l\le m$.

Estimators are constructed in several steps.

\medskip
\noindent
\textbf{Step 1. Estimators for the  parameters $(\beta,\underline{\mu},\overline{\mu},\underline{\sigma}^2)$:}
\quad

\begin{enumerate}
\item[(i)]
  For  each block  $1\leq l\leq m$  with data
  $B_l=\{(x_i,y_i)\}_{l\le i\le l+n-1}$, we run an ordinary LSE procedure using the
  standard regression model
  \[  y_i = \beta_l x_i +\mu_l + \varepsilon_{l,i}, \quad  l\le i \le l+n-1.
  \]
  Let $(\hat{\beta}_l,\hat{\mu}_l)$ be the obtained estimates for the
  regression parameter and mean parameter.

  Denote by $z_i=y_i-\hat{\beta}_l x_i-\hat{\mu}_l$ the corresponding
  residuals. Define the mean squared error from the $l$th block
  by
  $$
  \hat{\sigma}^2_l=\frac1{n-1}  \sum_{i=l}^{l+n-1} z_i^2 .
  $$

\item[(ii)]
  Find  the  block $\hat k$ with minimum mean squared error, that is,
  \[ 
     {\hat k} = \arg    \min_{1\leq l\leq m} \hat{\sigma}^2_l.
  \]
  Let
  \begin{align*}
    w_i & =y_i-\hat{\beta}_{\hat k}x_i,\quad 1\leq i\leq T, \\
    \tilde{\mu}_l & =\frac1n \sum_{i=l}^{l+n-1}w_i,\quad    1\leq l\leq m.
  \end{align*}

  We introduce the following  estimators.
  \begin{itemize}
  \item
    The lower  and upper means $\{\underline{\mu},
    \overline{\mu}\}$ are  estimated, respectively,  by
    \begin{equation}
      \hat{\underline{\mu}}=\min_{1\leq l\leq m} \tilde{\mu}_l,\qquad
      \hat{\overline{\mu}}=\max_{1\leq l\leq m} \tilde{\mu}_l.
      \label{eq:est-means}
    \end{equation}
  \item
      The regression parameter $\beta$ and the lower variance
      ${\underline\sigma}^2$ are  estimated by
      \begin{align}
        \hat \beta & = {\hat\beta}_{\hat k},        \label{eq:est-beta}        \\
        \underline{\hat\sigma}^2 & = {\hat\sigma}^2_{\hat k},  \label{eq:est-lower-variance}
    \end{align}
    that is, the regression estimators from the minimum mean squared  error block    $\hat k$.
  \end{itemize}

  Later, we will show that under appropriate conditions,
  the estimators
  $(\hat\beta,\hat{\underline{\mu}},\hat{\overline{\mu}},\underline{\hat\sigma}^2)$
  converge to  $(\beta,\underline{\mu},\overline{\mu},\underline{\sigma}^2)$ as $n\to \infty$ and $K\to \infty$.
\end{enumerate}

\noindent
\textbf{Step 2. Estimator for the upper variance $\overline{\sigma}^2$:}
\quad
 To estimate the upper variance $\overline{\sigma}^2$, we need to
 remove the mean uncertainty which is present in  the intermediate
 residuals $w_i=y_i-\hat{\beta}_{\hat k} x_i,\ 1\leq i\leq T$.
 Let $n_1<n$ be  a small window size and $P=T/n_1$ (in practice,
 values like $n_1=10,20,40$ are recommended).
 Let
 $$
 \tilde{w}_i=w_i-\frac1{n_1} \sum_{i=1+(j-1)n_1}^{jn_1}w_i,\quad   1+(j-1)n_1\le i \leq jn_1,\quad 1\leq j\leq P.
 $$
 This  steps centralizes the data over a local window, and is expected
 to remove the fluctuation (uncertainty) about observation means.
 Define, for $1\le l\le m$,
 $$
 \hat{\tilde{\sigma}}^2_l=\frac1{n-1}  \sum_{i=l}^{l+n-1} \tilde{w}_i^2.
 $$
 Finally we estimate the upper variance by
 \begin{equation}
   \label{eq:u-variance}
   \hat{\overline{\sigma}}^2 = \max_{1\leq j\leq m}\hat{\tilde{\sigma}}^2_j.
 \end{equation}

\medskip

The construction of the estimators above is motivated by  the
following observations.

\begin{enumerate}
\item[(i)\ ]
  When two groups of samples,  with respective sample
  means $(\mu_1,\mu_2)$ and sample variances
  $(\sigma^2_1,\sigma^2_2)$, are merged to  one group,  the mean of
  the resulting group takes value in the interval
  $[\mu_1\wedge \mu_2,    \mu_1\vee \mu_2]$; its variance is larger than $\sigma^2_1\wedge \sigma^2_2$.
\item[(ii)]
  If the two  groups have a  same sample mean $\mu$ and different sample
  variances $(\sigma^2_1,\sigma^2_2)$, the variance of the merged
  group belongs to the interval
  $[\sigma^2_1\wedge \sigma_2^2,\sigma^2_1\vee \sigma^2_2]$.
\end{enumerate}

Furthermore,  with reference to   the data generation  process \eqref{eq:DGP},
consider a  data group of length $n_0$,
$A_j=\{(x_i,y_i)\}_{i=1+(j-1)n_0}^{jn_0}$,
where $ 1\leq j\leq K$.  When $n\leq n_0$, there exists a moving
group  $B_l=\{(x_i,y_i)\}_{i=l}^{l+n-1}$, where  $1\leq l\leq m$,
such that $B_l\subset A_j$.
We use the ordinary LSE to estimate the regression and mean
parameters within
each of the data blocks  of $\{B_l\}_{l=1}^m$, and obtain $m$
estimates  for the regression coefficient $\beta$ and the
corresponding mean squared errors.
Based on observation (i), we can use the minimum mean squared error
from these $m$ data blocks as an estimator for the
minimum variance of the data groups $\{A_j\}_{j=1}^K$. This is done
with block $\hat k$ and the mean squared error $\hat\sigma_{\hat  k}^2$ from this block.
Then, by Theorem \ref{the:u-lse} and Corollary \ref{coro:1}, we can obtain the estimation $\hat{\beta}_{\hat{k}}$ for $\beta$ based on the block $\hat k$.

The next  question is  to estimate
$(\underline{\mu},\overline{\mu},\overline{\sigma})$ via the $m$ sets
of residuals $\{C_l\}_{l=1}^m$, where $C_l=\{w_i=y_i-\hat{\beta}_{\hat k} x_i\}_{i=l}^{l+n-1}$.
By observation (i), we can calculate  means $\{\tilde\mu_l\}$
of  these $m$ sets of
residuals, and their minimum and maximum values will be a good
estimator for
$\min\limits_{1\le j\le K}\eta_j$
and
$\max\limits_{1\le j\le K}\eta_j$, respectively. As the latter values
converge to the lower and upper mean,  $\underline\mu$ and
$\overline\mu$, respectively, when $K\to\infty$, these estimators are
consistent.

Finally for estimating the upper variance $\overline{\sigma}^2$ in \textbf{Step 2} ,
we first remove the mean uncertainty that is present  in the intermediate
residuals  $\{C_l\}_{l=1}^m$ by using local averaging over smaller blocks
of size $n_1<n$.
Then, by observation (ii),  we can estimate  $\overline{\sigma}$
with the maximum value of the mean squared errors from blocks
$\{C_l\}_{l=1}^m$ after removing mean uncertainty.

The theoretical consistency of these estimators are established in the
following theorem.
\begin{theorem}\label{the-consistency}
  Consider the data generation process~\eqref{eq:DGP}, and assume
  that as  $K\to \infty$,
  \[
  (\min_{1\leq j\leq K}\eta_j,\max_{1\leq j\leq K}\eta_j,\min_{1\leq
    j\leq K}\sigma^2_j,\max_{1\leq j\leq K}\sigma^2_j)
  \longrightarrow
  (\underline{\mu},\overline{\mu},\underline{\sigma}^2,\overline{\sigma}^2).
  \]
  Assume also $n_0\le n$. Then as  $K\wedge n\to\infty$,
  \begin{enumerate}
  \item[(i)]
    the lower  variance estimator is strongly consistent, that is,
    $\hat{\underline{\sigma}}^2 \to \lsigma$, with probability $1$;
  \item[(ii)]
    the   estimator $\hat \beta$ for the regression parameter
    is strongly consistent, that is,
    $\hat\beta \to \beta$, with probability $1$;

  \item[(iii)]
    the lower and upper mean estimators in \eqref{eq:est-means} are strongly consistent, that
    is, $(\hat{\underline{\mu}},\hat{\overline{\mu}})$ converge to
    $(\underline{\mu},\overline{\mu})$ with probability $1$;
  \item[(iv)]
    the upper variance estimator in \eqref{eq:u-variance} is strongly consistent, that is,
    $\hat{\overline{\sigma}}^2 \to \overline{\sigma}^2$ with probability $1$.
  \end{enumerate}
\end{theorem}
The proof of the theorem is  given in Appendix~\ref{proof-the-consistency}.

\section{Simulation experiments}
\label{sec:simul}


Simulations are conducted to check the finite-sample  performance of the
Robust-LSE estimators proposed in Section~\ref{ssec:estim}.
The design for the data generation process~\eqref{eq:DGP} is as
follows: for  $1\le j\le K$,    $1+n_0(j-1)\leq i\leq n_0j$,
\begin{itemize}
\item
  $\eta_j$ takes value in $[0,5]$ uniformly, $\sigma_j$ takes value in
  $[0.1,1]$ uniformly. Define
  $$
  ({\eta}_{\min},{\eta}_{\max})=(\min_{1\leq j\leq K}\eta_j,\max_{1\leq j\leq K}\eta_j),\quad
  ({\sigma}_{\min},{\sigma}_{\max})=(\min_{1\leq j\leq K}\sigma_j,\max_{1\leq j\leq K}\sigma_j).
  $$
\item
  $\varepsilon_i\in \EN(0,\sigma_j^2)$;
\item
  $ y_i=x_i+\eta_j+\varepsilon_i,   \quad (\beta=1).$
\end{itemize}

Consider the estimators
$(\hat{\beta},\hat{\underline{\mu}},\hat{\overline{\mu}},\hat{\underline{\sigma}},\hat{\overline{\sigma}})$
defined in  \textbf{Steps 1 and 2} in Section~\ref{ssec:estim}.
We take  $(n_0,n,n_1)=(200,150,20)$ and varying
$T\in\{400,800,1600,3200\}$ (or equivalently, $K=T/n_0\in\{2,4,8,16\}$).

For each combination of $(T,n_0,n,n_1)$, we generate 500 independent
replications of the data set
$\{\eta_j,\sigma_j\}_{j=1}^K$ and errors
$\{\varepsilon_i\}_{i=1}^T$.
The average values of the parameters
$({\eta}_{\min},{\eta}_{\max},{\sigma}_{\min},{\sigma}_{\max})$ over
the 500 replications are denoted as
$(\bar{\eta}_{\min},\bar{\eta}_{\max},\bar{\sigma}_{\min},\bar{\sigma}_{\max})$.

Table~\ref{table:1} reports empirical statistics for the Robust-LSE
estimators
and for comparison purpose, the ordinary LSE estimators.
For each case, we calculate the average and standard error for the two
estimators of  $\beta$.
The method Robust-LSE indeed provides  a better
estimator $\hat\beta$ than the ordinary LSE,
with smaller standard errors for $K\in\{2,4\}$ and comparable standard
errors for $K\in\{8,16\}$.
Note that, we have taken
the parameters $\{\eta_j,\sigma_j\}_{j=1}^K$ uniformly from some
intervals.
 The induced mean and variance uncertainties are less severe when the
  number of groups  $K$ grows because in this case, an averaging effect
  appears to reduce such uncertainties, and thus the ordinary LSE method
  is able to provide an accurate estimate for the regression parameter.
  However, if the uncertain mean and variance values $\{\eta_j,\sigma_j\}$
  do not obey any clearly defined distributions (as done here), the performance of the
  ordinary LSE is likely to worsen.
Furthermore by construction, the  Robust-LSE method provides
consistent estimators
$(\hat{\underline{\mu}},\hat{\overline{\mu}},\hat{\underline{\sigma}}^2,\hat{\overline{\sigma}}^2)$
for
the  mean and variance uncertainty parameters in  the samples.

\begin{table}[H]
\centering
\caption{Empirical statistics of the Robust-LSE estimators and the
  ordinary LSE estimator from 500 replications. Average and standard
  errors are reported for $\beta$.  Parameters are
  $\beta=1$,  $(n_0,n,n_1)=(200,150,20)$ and $T\in\{400,800,1600,3200\}.$}
\label{table:1}
\begin{tabular}{cccc}
\toprule
 &$\hat{\beta}$ & $(\hat{\underline{\mu}},\hat{\overline{\mu}})$ & $(\hat{\underline{\sigma}},\hat{\overline{\sigma}})$ \\
  \midrule
  \multicolumn{4}{|c|}{ $T=400,\ (\bar{\eta}_{\min},\bar{\eta}_{\max})=(1.7405,3.2692),~(\bar{\sigma}_{\min},\bar{\sigma}_{\max})
  =(0.3991,0.7055)$} \\
  Robust-LSE & 0.9729   &(1.7479,3.2705)    &(0.3742,0.7129) \\[-3mm]
             & (0.5151)   &&\\
  LSE        & 1.0299    &2.4562        &0.7352 \\[-3mm]
             & (1.4121)   &&\\
  \multicolumn{4}{|c|}{  $T=800,\ (\bar{\eta}_{\min},\bar{\eta}_{\max})=(0.9819,3.9324),
    ~(\bar{\sigma}_{\min},\bar{\sigma}_{\max})  =(0.2890,0.8275)$} \\
  Robust-LSE &  0.9820    & (0.8607,3.9908)    &(0.2703,0.8399) \\[-3mm]
             & (0.3583)    &&\\
  LSE      & 0.9839      &2.6162        &1.1562  \\[-3mm]
           & (0.5803)    &&\\
  \multicolumn{4}{|c|}{ $T=1600,\ (\bar{\eta}_{\min},\bar{\eta}_{\max})=(0.5733, 4.4535),~(\bar{\sigma}_{\min},\bar{\sigma}_{\max})  =(0.1975,0.8956)$} \\
  Robust-LSE &   0.9994    &(0.4028, 4.5665)    &(0.1861,0.9164) \\[-3mm]
             &  (0.2387)   &&\\
    LSE     & 1.0147      &2.5216      &1.3733 \\[-3mm]
         &  (0.2105)     &&\\
  \multicolumn{4}{|c|}{ $T=3200,\ (\bar{\eta}_{\min},\bar{\eta}_{\max})=(0.2987,4.6994),~(\bar{\sigma}_{\min},\bar{\sigma}_{\max})
  =(0.1525, 0.9491)$} \\
  Robust-LSE & 0.9916    &(0.0730,4.9132)    &(0.1427,0.9776) \\[-3mm]
  &  (0.1153)  &&\\
  LSE     &  0.9948   &    2.4905              &1.4729 \\[-3mm]
  &   (0.0812) &&\\
   \bottomrule
  \hline
\end{tabular}
\end{table}

\begin{figure}[H]
  \centering
    \caption{Samples of  regression lines from the ordinary LSE  and
      from the  minimum mean squared error block $\hat{k}$ LSE.\label{figs:com}}
\includegraphics[width=2.95 in]{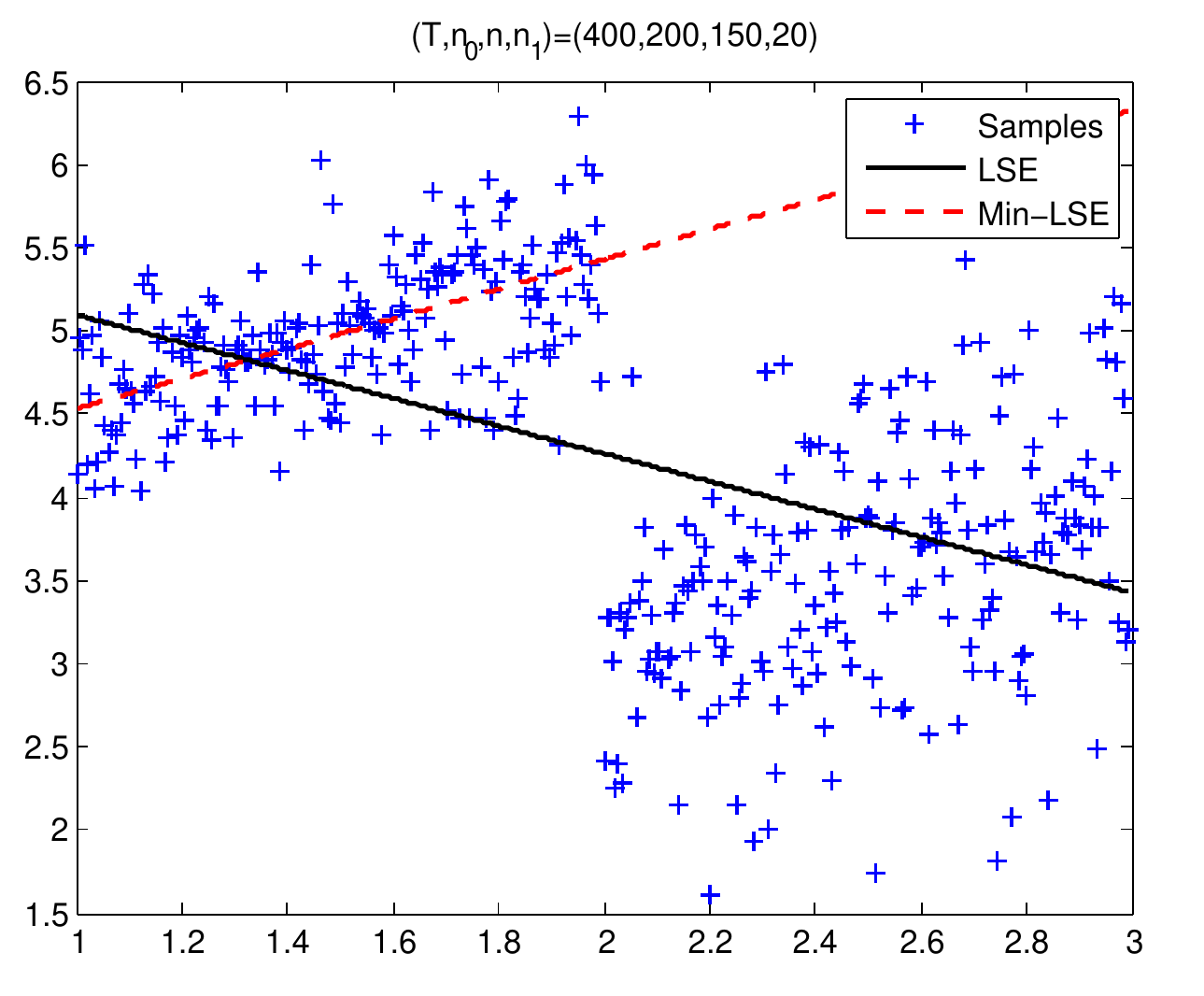}
\includegraphics[width=2.95 in]{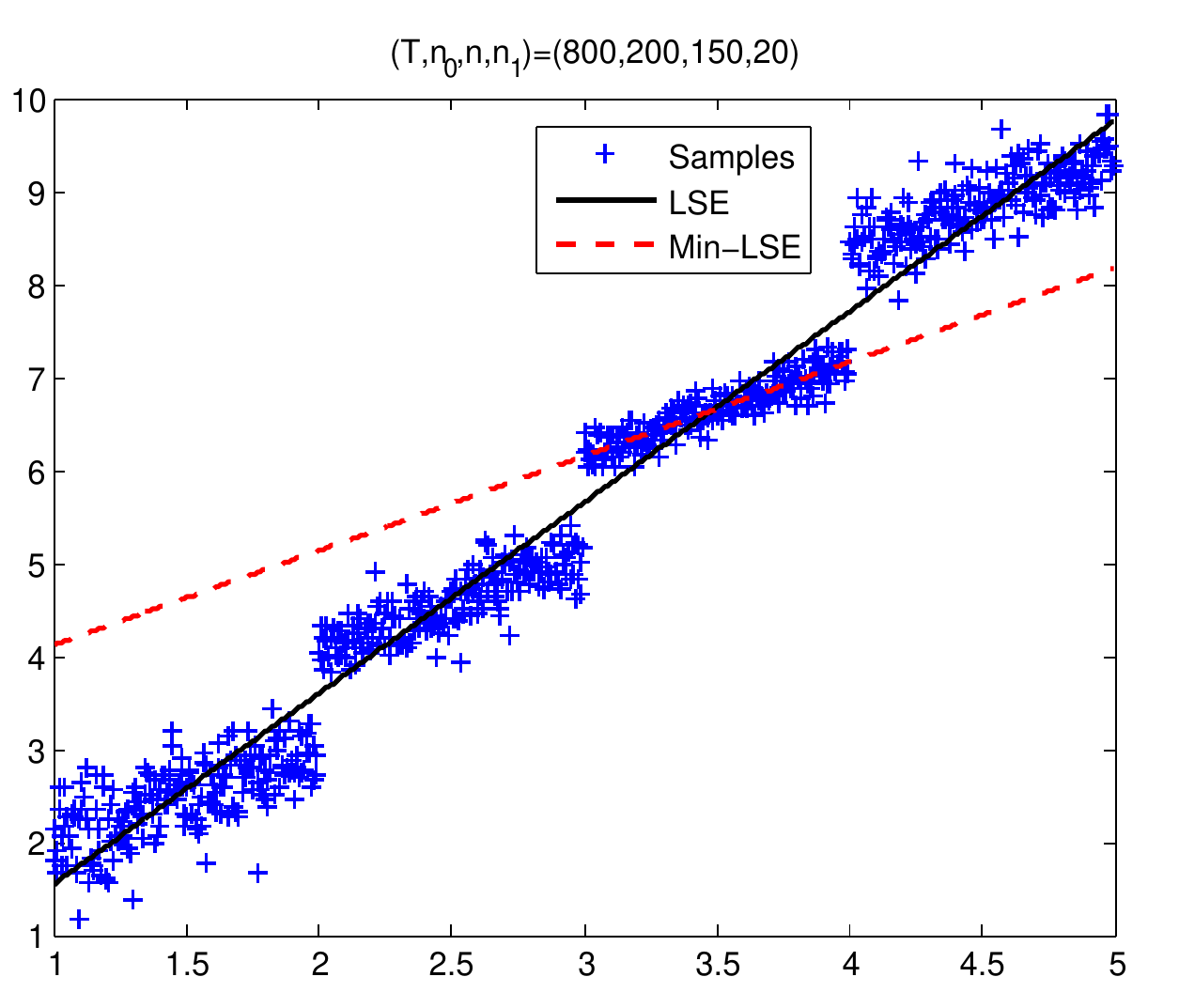}
\includegraphics[width=2.95 in]{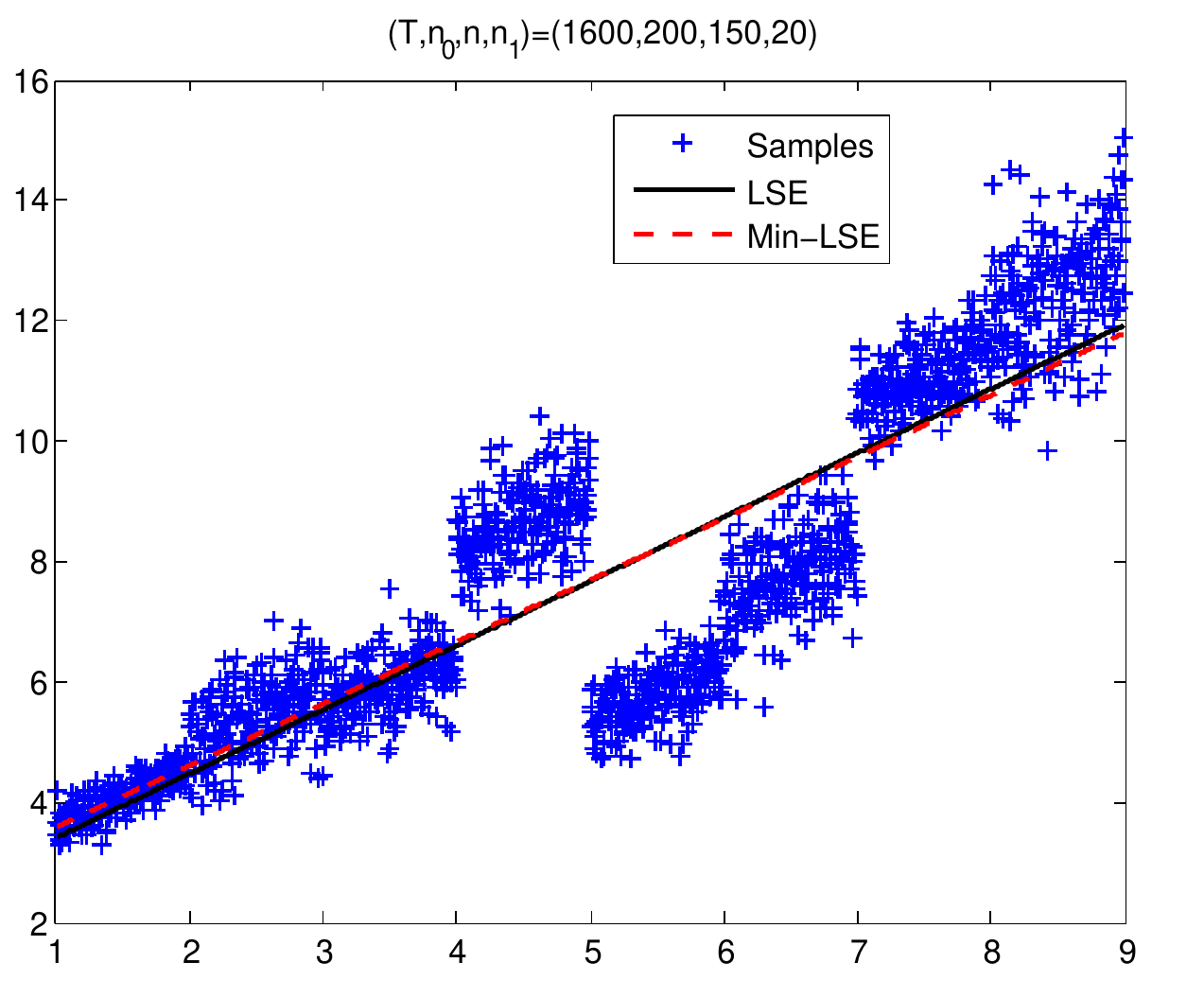}
\includegraphics[width=2.95 in]{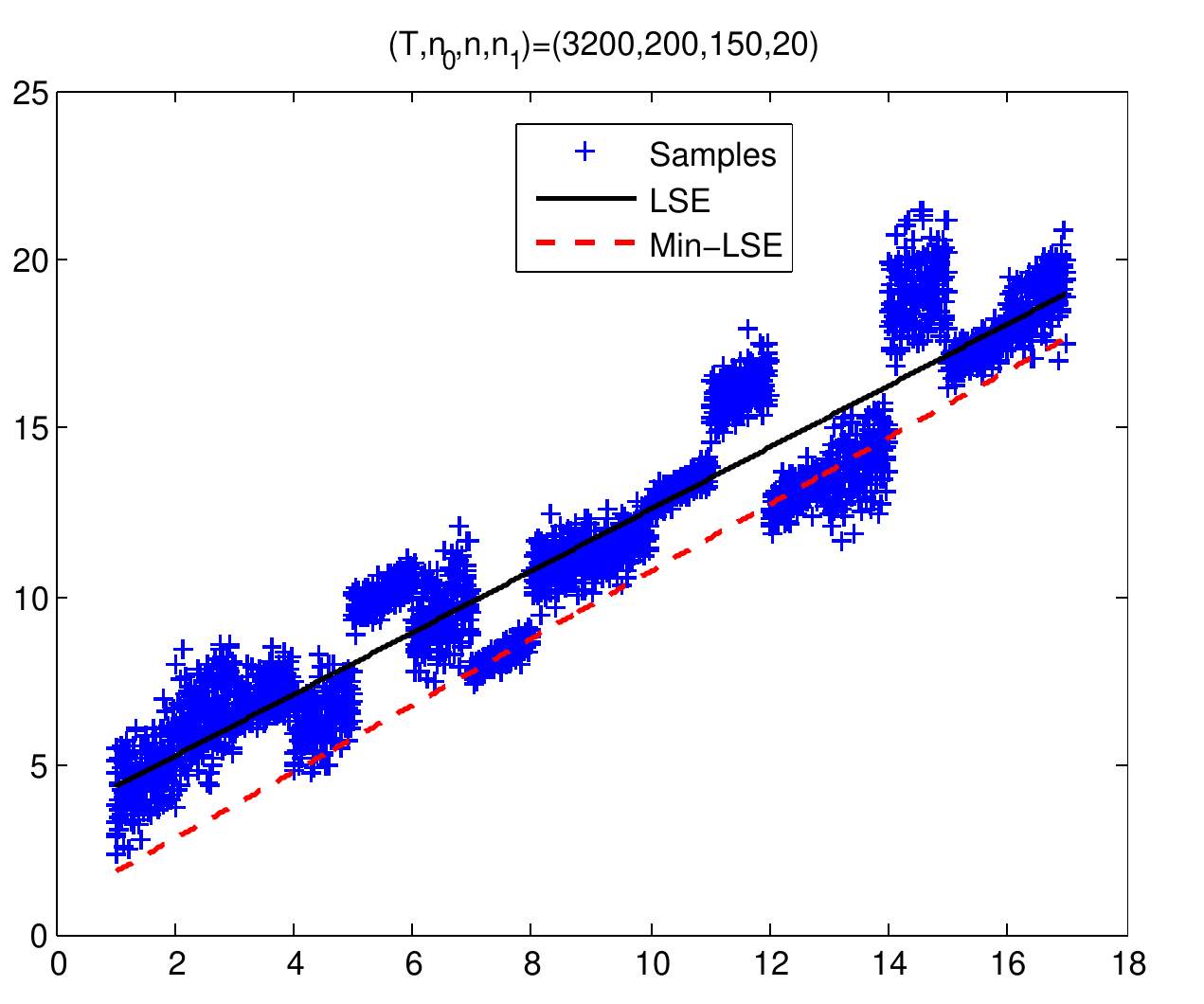}
\end{figure}

In Figure \ref{figs:com}, we plot samples of regression  lines for the
ordinary LSE and from the minimum mean square error  block $\hat{k}$ LSE with
parameters $(\hat\beta_{\hat k},{\hat\mu}_{\hat k})$ given in the
Robust-LSE method, respectively, and for sample size
$T\in\{400,800,1600,3200\}$. It is clear  that lines  from
minimum mean square error  block $\hat{k}$ LSE focus  on the
sub-samples with minimum mean squared errors, while lines from the ordinary LSE
focus  on the whole sample.
This explains why in general the Robust-LSE method can provide a
better estimator for the regression parameter ${\beta}$.

Now, we fix the value of $(T,n_0,n_1)=(1600,200,20)$, and verify
properties  of the Robust-LSE estimators when the
block length $n$ grows.
We take $n=60,80,160,200$.
Table \ref{table:2}
reports empirical statistics from 500 replications. From the estimators
$(\hat{\beta},\hat{\underline{\mu}},\hat{\overline{\mu}},\hat{\underline{\sigma}},\hat{\overline{\sigma}})$,
it is again observed that the method Robust-LSE performs better than
the  ordinary LSE.
Furthermore, the convergence of the Robust-LSE estimators
$(\hat{\beta},\hat{\underline{\mu}},\hat{\overline{\mu}},\hat{\underline{\sigma}},\hat{\overline{\sigma}})$
is verified when $n$ grows.
In Figure~\ref{figs:com1}, we show as in Figure~\ref{figs:com},
sample regression lines from the ordinary LSE and the
minimum mean squared error  block $\hat k$ from the Robust-LSE
method.  We observe that the latter can catch the groups with
minimum variance in all the tested cases of block length $n$.

\begin{table}[H]
\centering
\caption{
  Empirical statistics of the Robust-LSE estimators and the
  ordinary LSE estimator from 500 replications. Average and standard
  errors are reported for $\beta$.  Parameters are
  $\beta=1$,  $(T,n_0,n_1)=(1600,200,20)$ and $n\in\{60,80,160,200\}.$
  \label{table:2}}
\begin{tabular}{cccc}
\toprule
& $\hat{\beta}$ & $(\hat{\underline{\mu}},\hat{\overline{\mu}})$ & $(\hat{\underline{\sigma}},\hat{\overline{\sigma}})$ \\
  \midrule
     \multicolumn{4}{|c|}{
       $n=60,\ (\bar{\eta}_{\min},\bar{\eta}_{\max})=(0.5642,4.4594),~(\bar{\sigma}_{\min},\bar{\sigma}_{\max}) =(0.1965,0.8901)$} \\
    Robust-LSE         & 0.9877    &(0.1597,4.9009)    &(0.1633,0.9994) \\[-3mm]
             & (0.3214) &&\\
    LSE               & 1.0142    &2.4645        &1.3726 \\[-3mm]
    &  (0.2175) &&\\
    \multicolumn{4}{|c|}{ $n=80,\ (\bar{\eta}_{\min},\bar{\eta}_{\max})=(0.5683,4.4583),~(\bar{\sigma}_{\min},\bar{\sigma}_{\max})
  =(0.2013,0.9046)$} \\
    Robust-LSE        & 0.9950    &(0.3790,4.8098)    &(0.1755,0.9810)\\[-3mm]
        &  (0.2193) &&\\
    LSE                &1.0054   &2.4732    &1.3782\\[-3mm]
    &  (0.2168) &&\\
    \multicolumn{4}{|c|}{ $n=160,\ (\bar{\eta}_{\min},\bar{\eta}_{\max})=(0.5857,4.4779),~(\bar{\sigma}_{\min},\bar{\sigma}_{\max})
  =(0.2034,0.8995)$} \\
    Robust-LSE        & 1.0028    &(0.4650,4.5696)   &( 0.1926,0.9123) \\[-3mm]
        &  (0.1769) &&\\
    LSE               & 0.9893  &2.5715       &1.3637 \\[-3mm]
    & (0.2126)&&\\
    \multicolumn{4}{|c|}{ $n=200,\ (\bar{\eta}_{\min},\bar{\eta}_{\max})=(0.5610, 4.4350),~(\bar{\sigma}_{\min},\bar{\sigma}_{\max})
  =(0.2024,0.8955)$} \\
    Robust-LSE        & 1.0016    &(0.5519,4.5037)    &(0.1959,0.8887)\\[-3mm]
    & (0.1039)&&\\
    LSE               & 0.9980    &2.5469        &1.3712  \\[-3mm]
               &  (0.2175)&&\\
   \bottomrule
  \hline
\end{tabular}
\end{table}

\begin{figure}[H]
  \centering
  \caption{ Samples of regression lines
    from  the ordinary LSE, and from the
    minimum mean squared error block $\hat{k}$ LSE.
    Parameters are
    $\beta=1$,  $(T,n_0,n_1)=(1600,200,20)$ and $n\in\{60,80,160,200\}.$
    \label{figs:com1}}
\includegraphics[width=2.95 in]{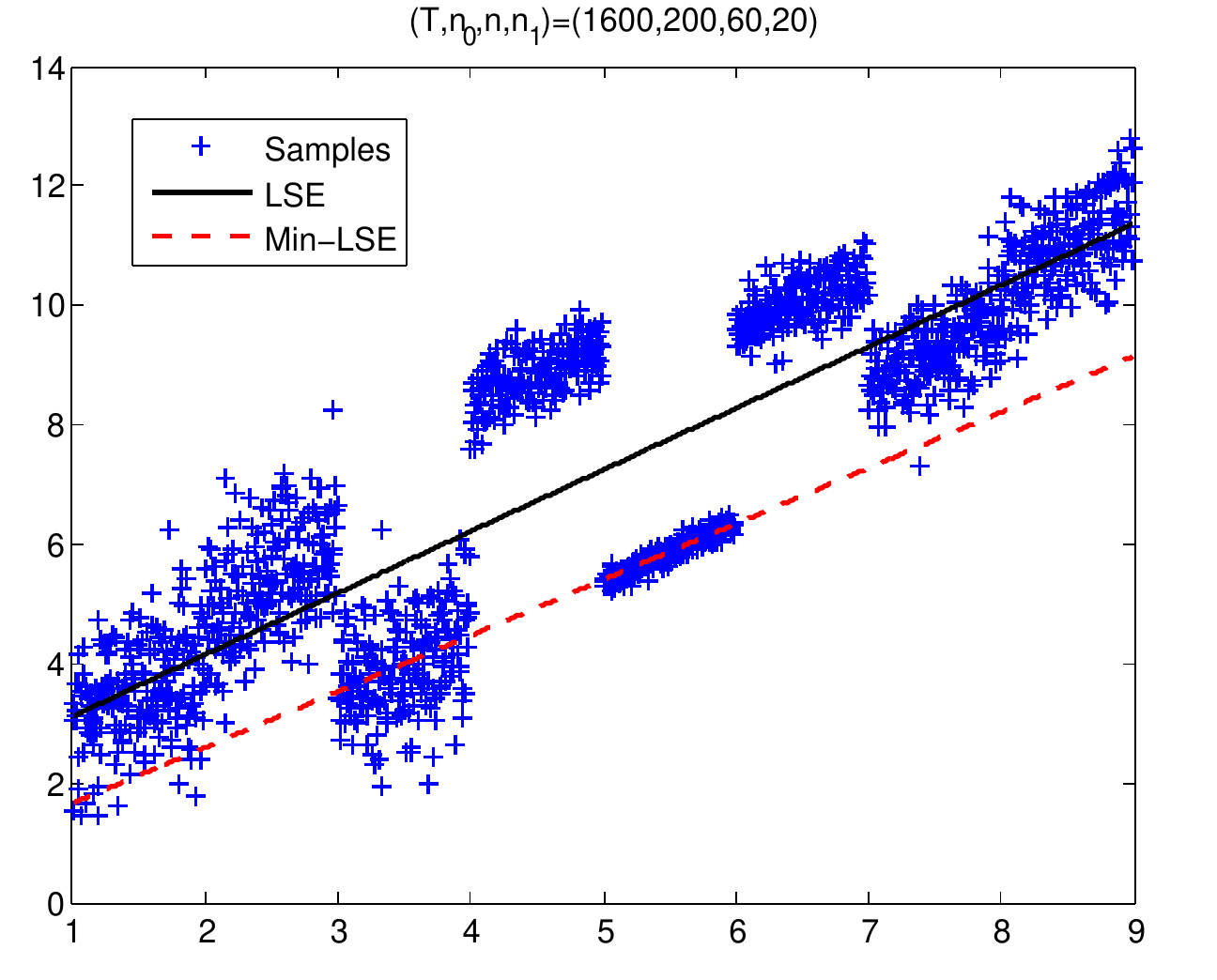}
\includegraphics[width=2.95 in]{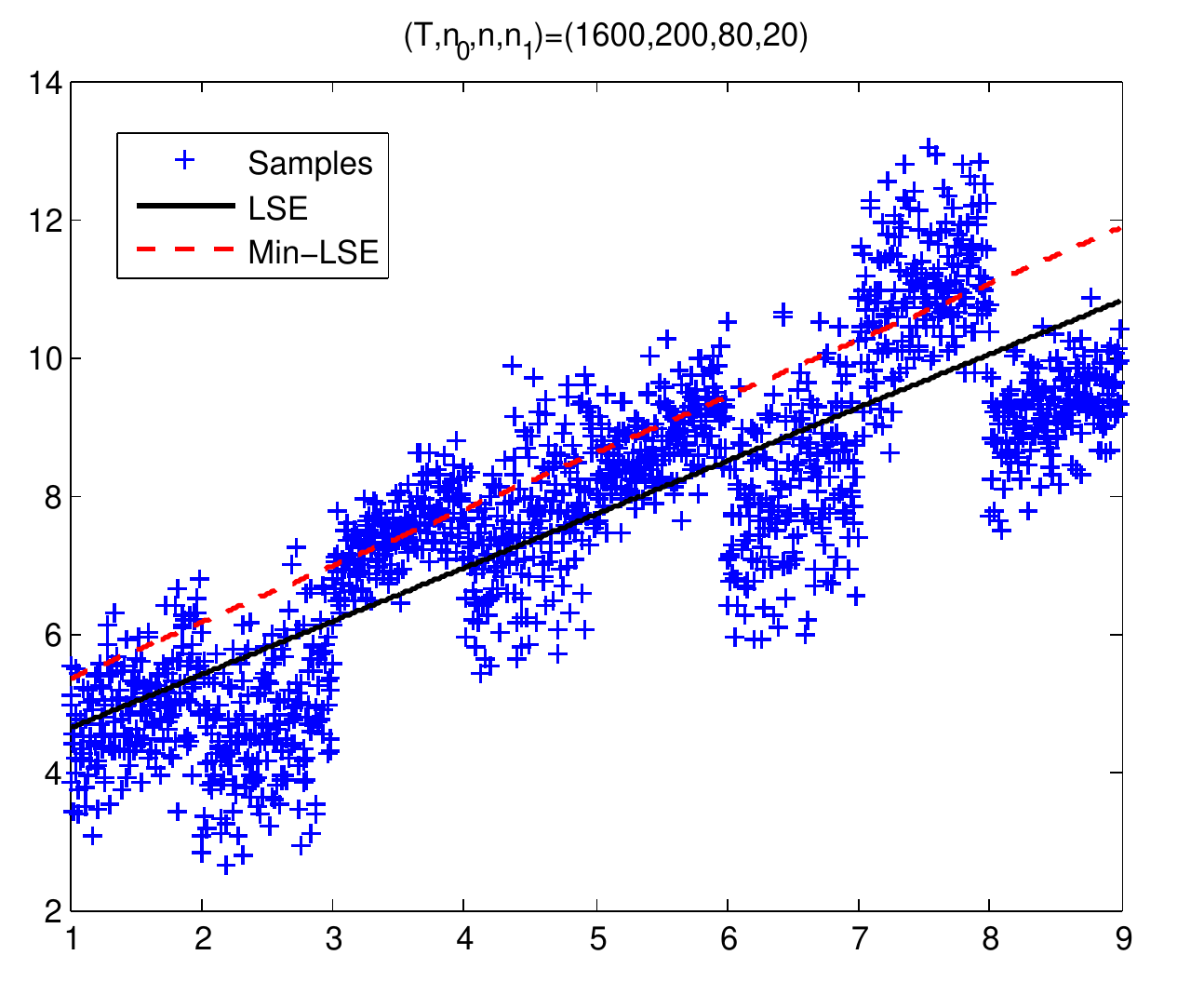}
\includegraphics[width=2.95 in]{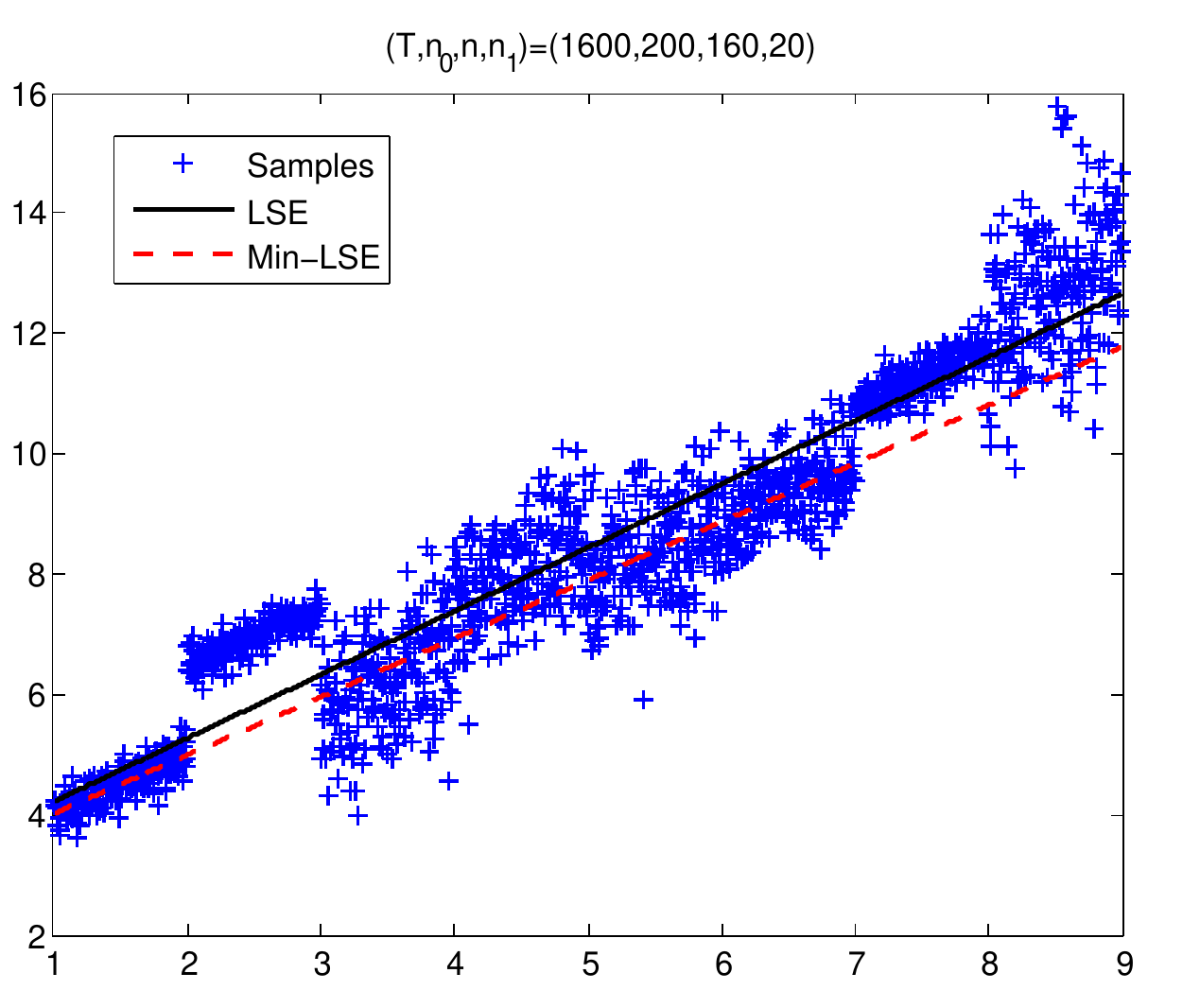}
\includegraphics[width=2.95 in]{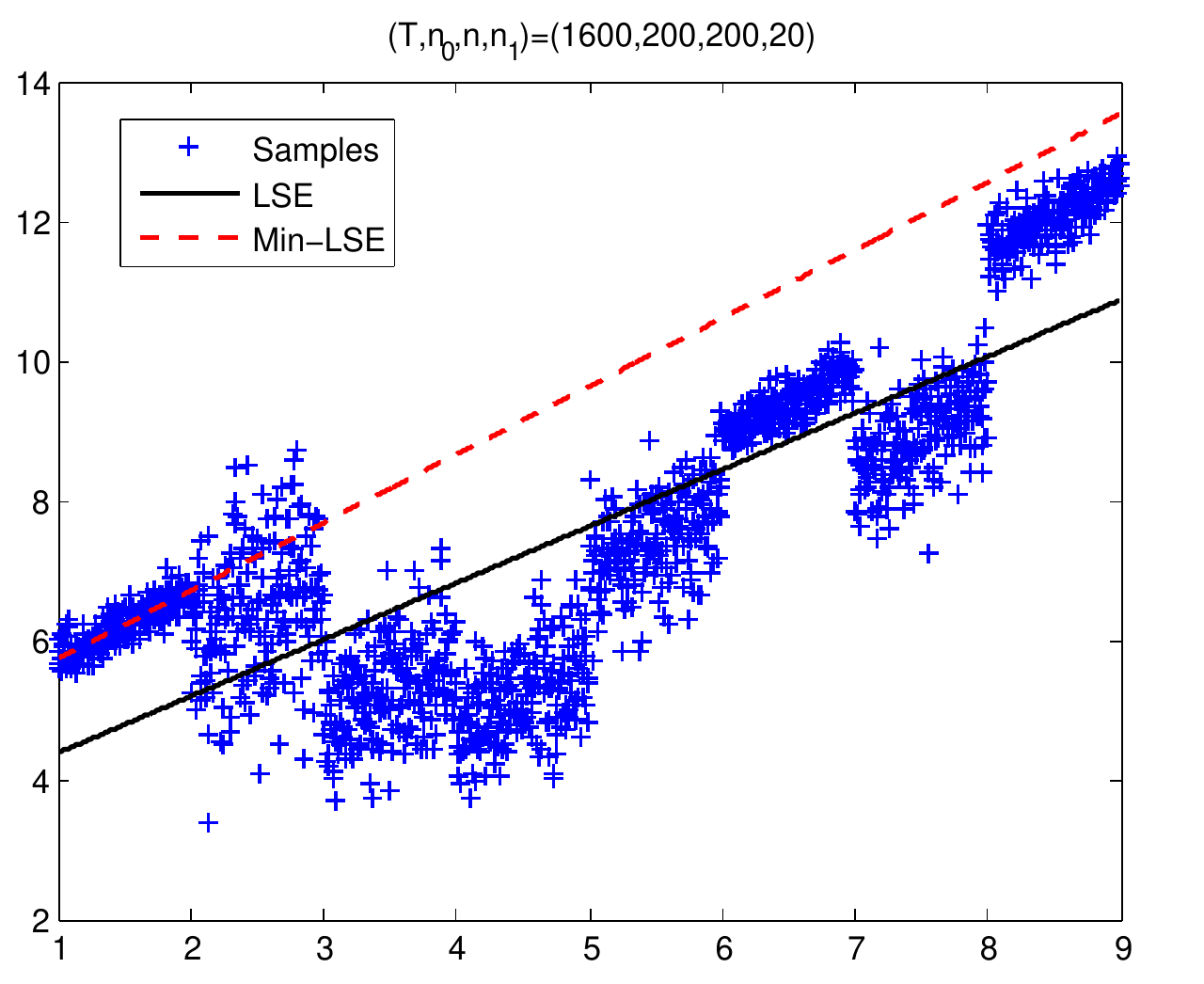}
\end{figure}

\section{Applications}\label{sec:app}
\subsection{Robust regression}

We apply the Robust-LSE estimators to the traditional robust
regression problem. Precisely, we compare our method with 
a benchmark robust regression estimator, namely  the MM estimator.
Actually, \cite{Yu17} has given an extensive review and comparison of
the existing robust regression estimators under various scenarios of
model contamination. Overall two estimators perform better than the
other competitors, namely the MM estimator \citep{Yohai87} and the
REWLSE estimator  \citep{Gervini02}. Since  these two best performers
are close each other, we chose the MM estimator as a reference in this
study.

Following a classical setting in the literature on robust regression,
we  consider a simple  linear model with contamination of the form
$Y=X+\varepsilon$, with samples
$\{(x_i,y_i)\}_{i=1}^T$ where
$x_i=1+0.01*i,\ 1\leq i\leq T$,  and 6 scenarios for
the errors $\{\varepsilon_i\}$: for $1\le m\le 6$,

\begin{description}
\item{Scenario $m$:} \quad $\varepsilon_i\in N(0,1),\ 1\leq i\leq a_m*T$, $\varepsilon_i\in N(0,100), a_m*T< i\leq T$.
  \\ Here $a_m\in\{0.95,0.90,0.80,0.85,0.70,0.60,0.50\}$, and
$1-a_m$ is  referred as the contamination rate of the base
  standard normal errors by a normal error with larger variance 100.
\end{description}

Under each scenario,
we generate 500 replications of the data,
and calculate the ordinary LSE, the MM estimator and the Robust-LSE
estimator for the regression parameter $\beta$.
Table~\ref{table:MSE} reports the MSEs of the estimators from 500
replications; a companion plot for  these MSEs is given  at the bottom of the
table.
We can see that in general, the ordinary LSE has a large MSE.
The Robust-LSE and the MM estimators 
have  almost identical performances for scenarios  1 and  2 with light
contamination. In contrast  for scenarios 3, 4 and 5 with heavier
contamination, the Robust-LSE clearly outperforms the MM estimator:
especially in the last case with 50\% contamination, the MM estimator
shows a breakdown with a MSE almost the double of the one from the
ordinary LSE (about 10 times of the one from the Robust-LSE
estimator).

\begin{table}[ptb]
\centering
\caption{Empirical MSEs of the ordinary LSE, MM and Robust-LSE
  estimators for the regression coefficient $\beta$. Sample size  $T=200$ with 500 replications.}
\label{table:MSE}
\begin{tabular}{cccc}
\toprule
   Scenario     & MM & LSE   & Robust-LSE  \\
  \midrule
   1 &0.0205   &0.2075    &0.0228 \\
   2 &0.0264   &0.3459    &0.0237\\
   3 &0.0603   &0.5733    &0.0405\\
   4 &0.1823   &0.6875    &0.0547\\
   5 &0.6079   &0.6970    &0.1060\\
   6 &1.6230   &0.7033    &0.1667\\
   \bottomrule
  \hline
\end{tabular}
\\
  \includegraphics[width=3.8 in]{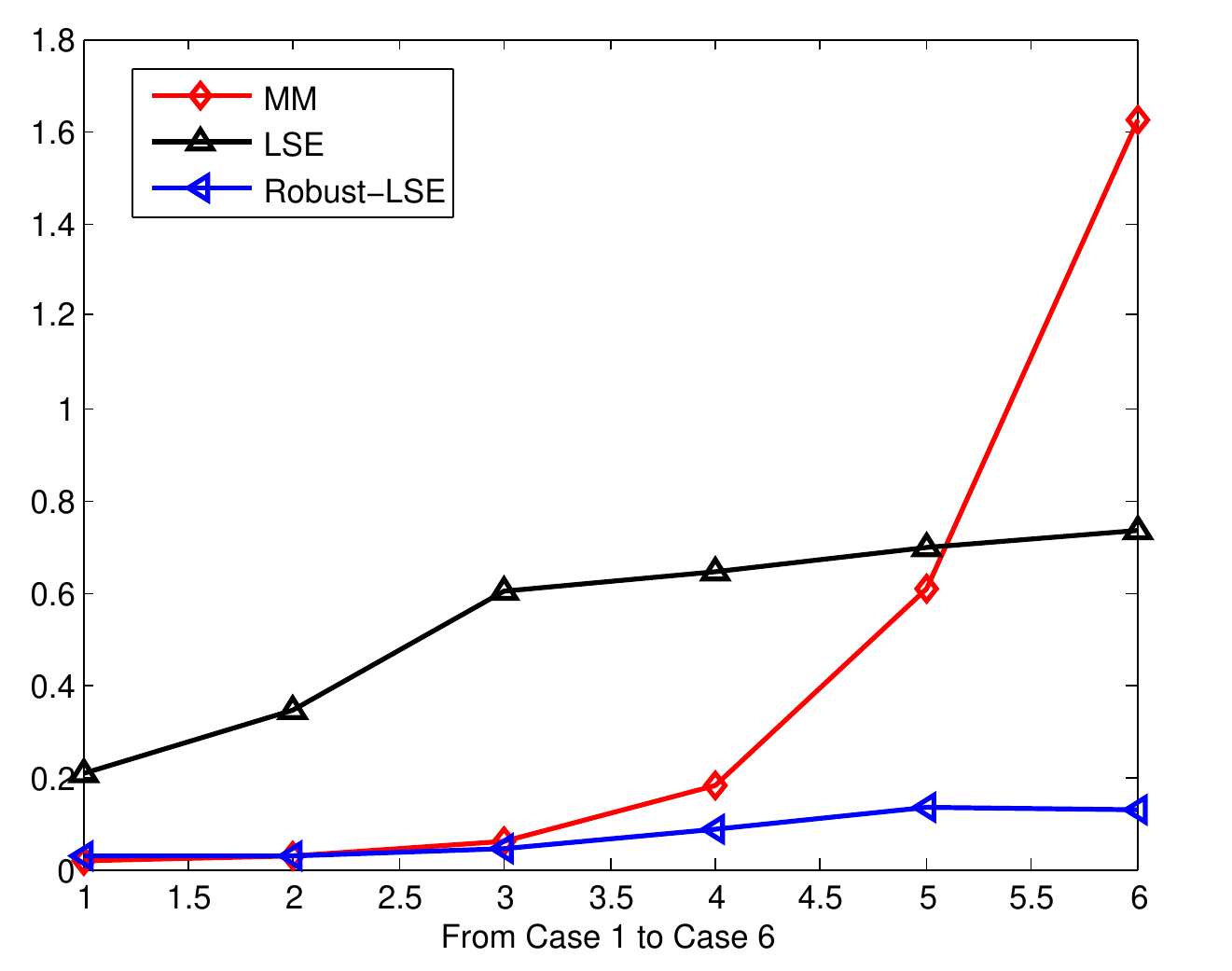}
\end{table}

Next we examine the  large sample behaviour of the three estimators by
gradually increasing the sample size from $T=200$ to $T=1000$.
Among the 6 scenarios of contamination, we report the results for
scenarios 1 and 4. The empirical MSEs are reported in
Table~\ref{table:MSE-large-T}, and displayed in a plot at its bottom.
We can see that the Robust-LSE estimator performs better than the
MM estimator in scenario 4 (medium contamination) while they are
similar under scenario 1 (light contamination) while being both
preferable than the ordinary LSE estimator. Besides, all the three
estimators show consistency when the sample size increases.

\begin{table}[H]
\centering
\caption{Empirical MSEs of the ordinary LSE, MM and Robust-LSE
  estimators for the regression coefficient $\beta$. Sample size
  $T\in \{200,400,600,800,1000\}$
  under scenarios 1 and 4  with 500 replications. }
\label{table:MSE-large-T}
\begin{tabular}{crccccc}
\toprule
 &   &$T=200$ & $T=400$   & $T=600$ & $T=800$ & $T=1000$ \\
  \midrule
  Scenario 1 & MM  &0.0205   &0.0024    &0.0007  & 0.0003  & 0.0002 \\
  & LSE    &0.2075  &0.0282  &0.0071  &0.0031  &0.0017 \\
  & R-LSE  &0.0228  &0.0025  &0.0007  &0.0003  &0.0001  \\
  Scenario 4 & MM &0.1823   &0.0230    &0.0064  &0.0029   &0.0013  \\
  &  LSE  &0.6875   &0.0833  &0.0251  &0.0116  &0.0052\\
  & R-LSE  &0.0547   &0.0059  &0.0020  &0.0007  &0.0004 \\
   \bottomrule
  \hline
\end{tabular}
\\

\includegraphics[width=3.0 in]{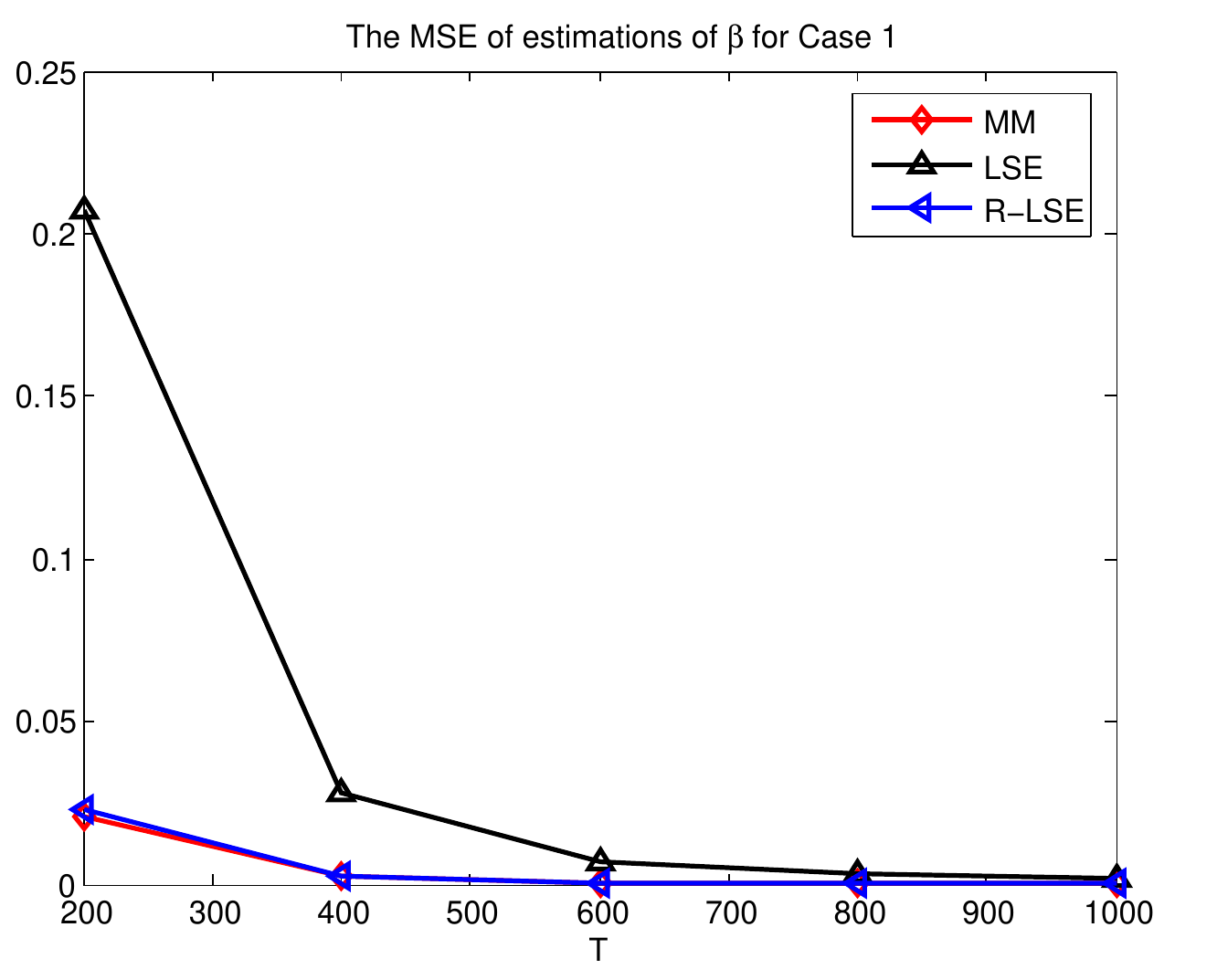}
\includegraphics[width=2.9 in]{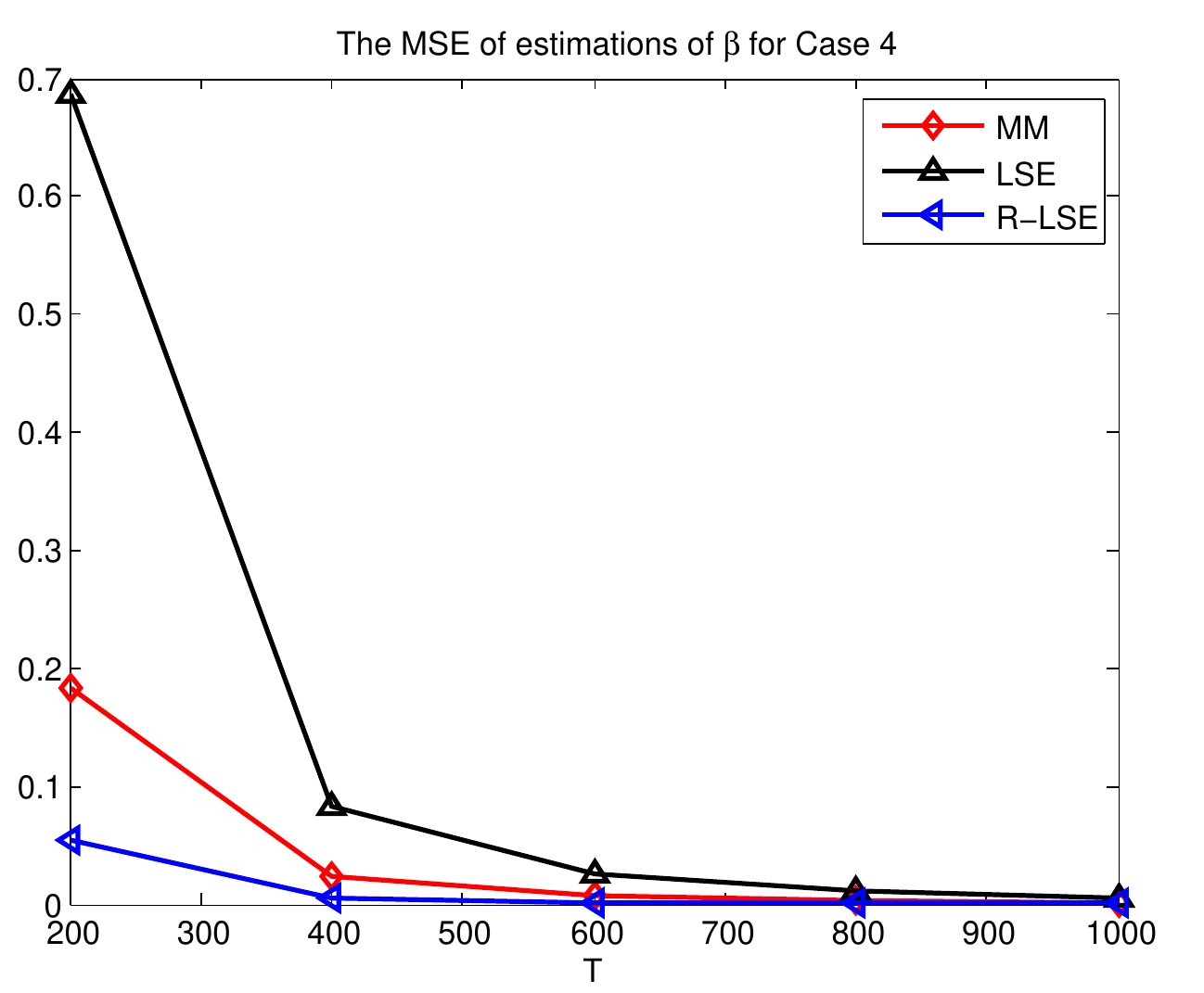}
\end{table}

\subsection{Regression under heteroscedastic errors}

In this section, we consider a special regression model under  heteroscedastic  errors:
$$
Y_{ij}=\beta X_{ij}+\varepsilon_{ij},\ 1\leq j\leq n_0,\ 1\leq i\leq K,\ T=n_0K,
$$
where $\varepsilon_{ij}\in \EN(0,\sigma_i^2),\ 1\leq j\leq n_0$. We set $\beta=1,K=10$, $X_{ij}=1+0.005(j+(i-1)n_0)$ and
$$
\{\sigma_i\}_{i=1}^{10}= \{ 0.6995,0.5851,0.3481,0.1304,0.7165,0.3344,0.4721,0.5211,0.1955,0.4851\}
$$
with $(\min_{1\leq i\leq 10}\sigma_i,\max_{1\leq i\leq 10}\sigma_i)=(0.1304,0.7165)$.
 This list of variances is quite arbitrary; their exact
  values have no particular meaning in our discussion.

The particularity here is that the model has only variance
uncertainty. We apply our Robust-LSE method,  without prior knowledge
about the heteroscedasticity of the data set, to obtain an estimation
for the regression parameter $\beta$ and the underlying minimum and
maximum volatility
$({\overline\sigma}_{\min},{\overline\sigma}_{\max})$.
Table~\ref{table:hetero} reports empirical averages of these estimates
from 500 replications. The corresponding ordinary LSE estimates are
also given for comparison.
Figure~\ref{figs:hetero} plots these empirical values.
We find that the Robust-LSE can provide an  estimator  for $\beta$
which is as good as the ordinary LSE; it can  also  provide
accurate estimations for the minimum and maximum volatilities while
the ordinary LSE cannot.

\begin{table}[H]
\centering
\caption{Heteroscedastic regression models with
  $(\beta,\bar{\sigma}_{\min},\bar{\sigma}_{\max})=(1,0.1304,0.7165)$. Averages
  of estimators from 500 replications and sample size $T\in\{500,1000,1500,2000\}.$}
\label{table:hetero}
\begin{tabular}{ccccc}
\toprule
& Parameters & $\beta$ & Min. volatility   & Max. volatility \\
  \midrule
$T=500$ &         R-LSE       &0.9770    &0.1213    &0.7844  \\
      &  LSE         &0.9981    &0.4845    &0.4845   \\
$T=1000$ &        R-LSE       &1.0066    &0.1246    &0.7691  \\
     &   LSE         &1.0001    &0.4864    &0.4864  \\
$T=1500$ &        R-LSE       &1.0050    &0.1258    &0.7544  \\
     &   LSE         &1.0002    &0.4858    &0.4858  \\
$T=2000$ &        R-LSE       &1.0013    &0.1267    &0.7466  \\
     &   LSE         &1.0000    &0.4854    &0.4854  \\
   \bottomrule
  \hline
\end{tabular}
\end{table}
 \begin{figure}[H]
  \centering
  \caption{Plots of empirical averages of the estimates in
    Table~\ref{table:hetero} (see captions there).\label{figs:hetero}}
\includegraphics[width=4.0 in]{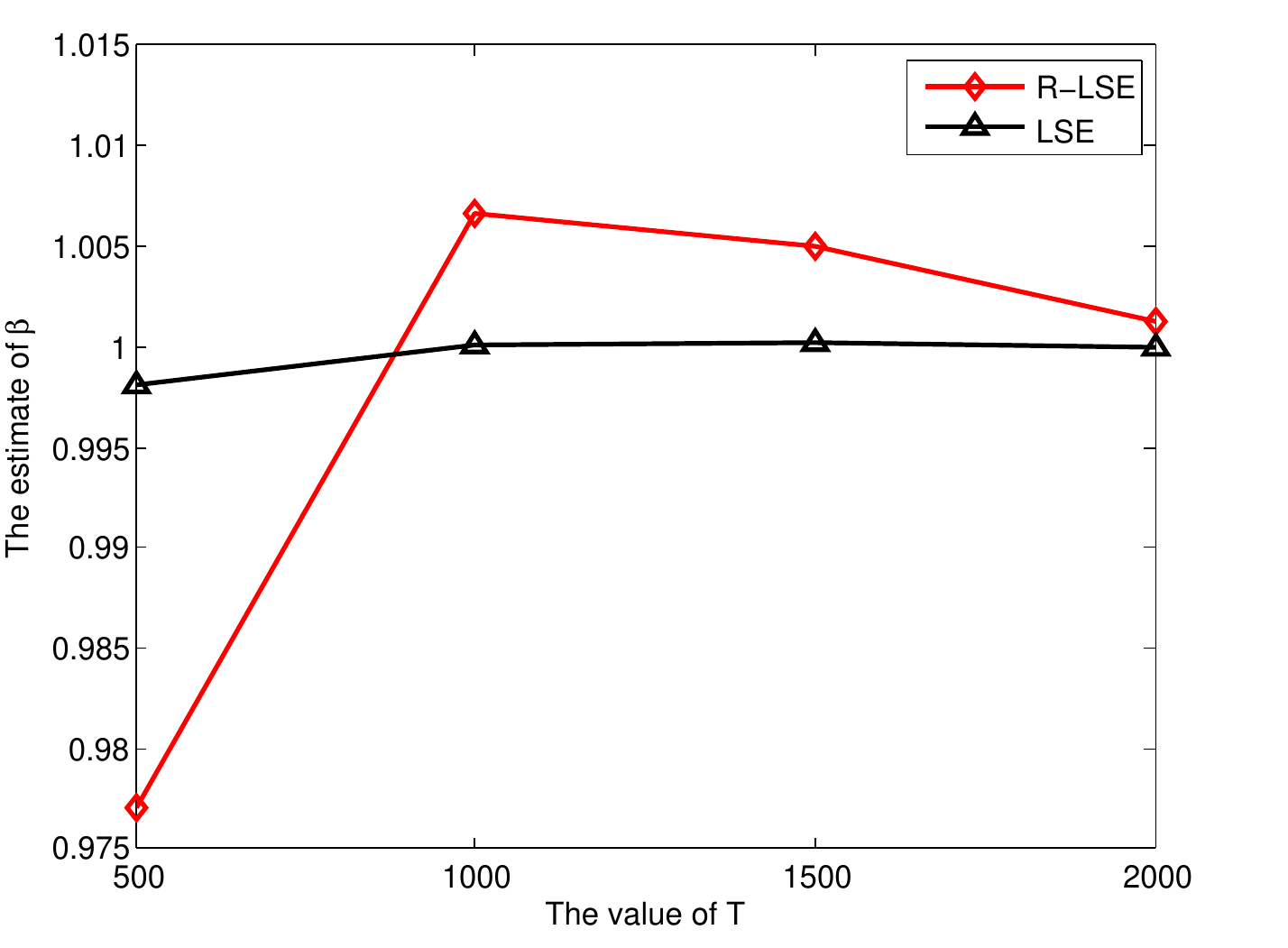}
\includegraphics[width=2.9 in]{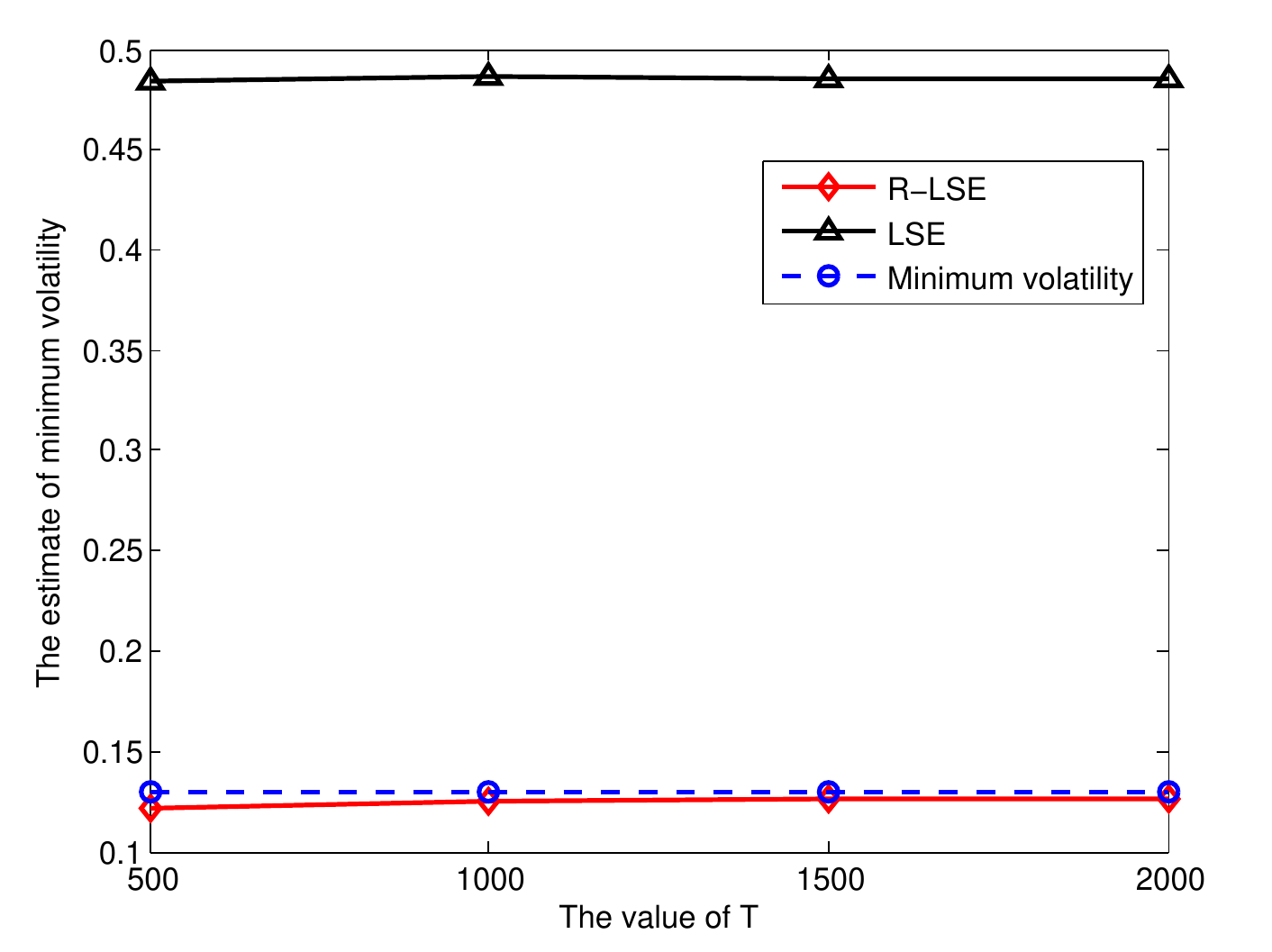}
\includegraphics[width=2.9 in]{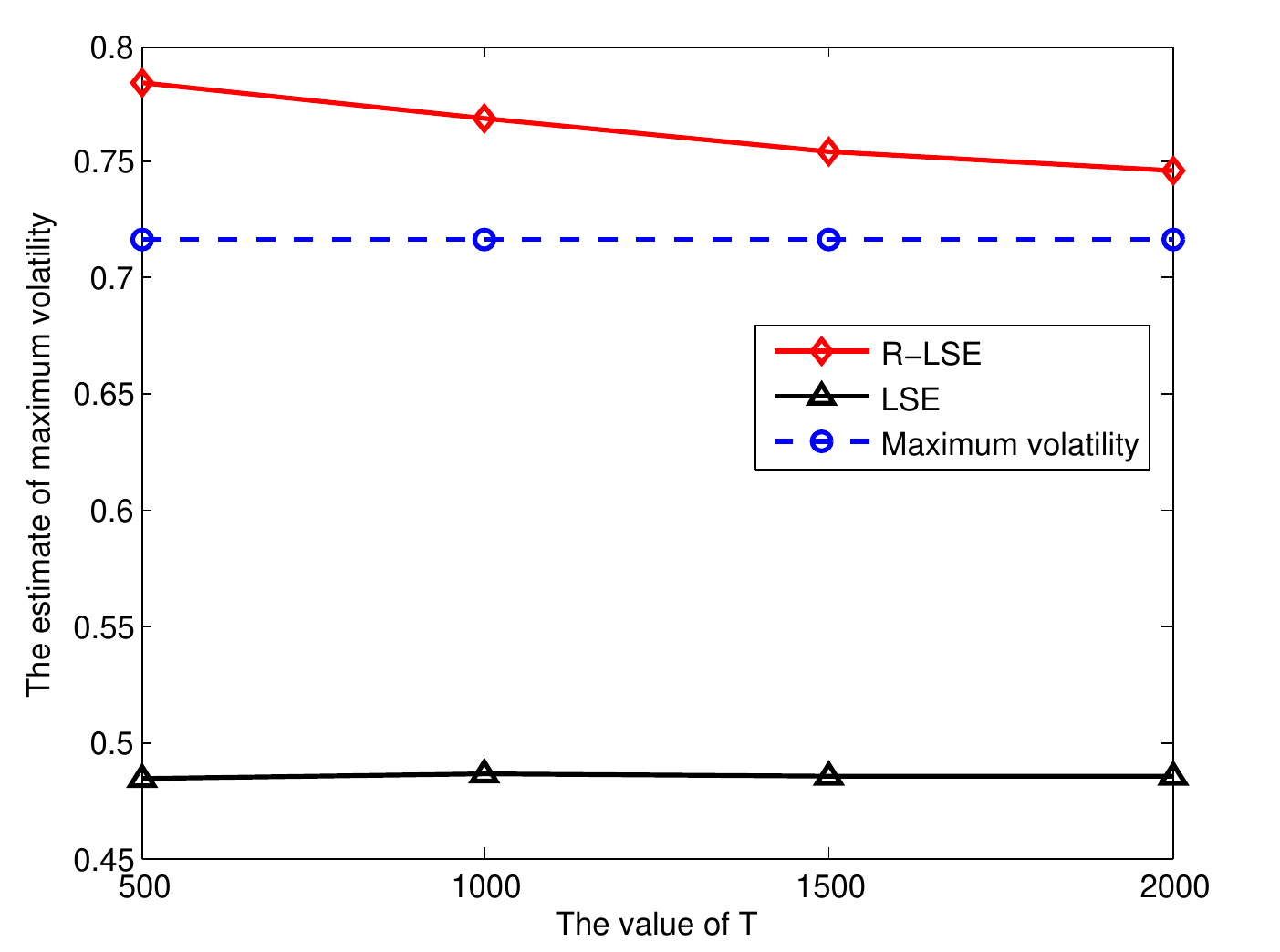}
\end{figure}

%
\subsection{Real data analysis}
\label{sec:real-data}
We consider a simple linear model:
\begin{equation}\label{eq:m1}
Y=\beta_1 X_1+\beta_2 X_2+\varepsilon,\ X_1,X_2\in \mathbb{R},
\end{equation}
where $\varepsilon$ satisfies a normal distribution $\EN(0,\sigma^2)$. In real market, it is important to select the factors for the linear regression model. However, we may not observe the factor $X_2$ and ignore it. Thus, it is possible that we consider the following model:
\begin{equation}\label{eq:m2}
Y=\beta_1 X_1+L+\varepsilon,\ X_1\in \mathbb{R},
\end{equation}
where $L$ is a constant. Note that, we can use the ordinal LSE to
obtain the coefficient of model (\ref{eq:m2}). Based on the
distribution-uncertain regression model (\ref{eq:line}), we use a mean uncertain term to represent the unknown factor $\beta_2 X_2$. The new model is
\begin{equation}\label{eq:m3}
Y=\beta_1 X_1+\eta+\varepsilon,\ X_1\in \mathbb{R},
\end{equation}
where $\eta$ takes value in a interval under sublinear expectation,
and $\varepsilon$ has the  $\EN(0,\sigma^2)$ distribution.

We analyze the  S\&P500 Index  to assess the performance of the models
(\ref{eq:m2}) and (\ref{eq:m3}). The daily closing price data of the
index covers the period  from Jan. 3, 2000 to July 17, 2020. We consider a first order autoregression version of models (\ref{eq:m2}) and (\ref{eq:m3}):
$$
X_{t+1}=\beta X_t+L+\varepsilon_{t+1},\quad X_{t+1}=\beta X_t+\eta+\varepsilon_{t+1}.
$$
\begin{table}[H]
\centering
\caption{Regression results of LSE and Robust-LSE under criterion $F_{0.01}(2,247)=4.6921$}
\label{table:SP500}
\begin{tabular}{clccc}
\toprule
 Year & Method  & $\beta$ & $R^2$   & $F$-statistic  \\
  \midrule
    201907--202007  &  R-LSE    &-0.4196    &0.1802   &\textbf{27.1386}  \\
                    &LSE        &-0.3592    &0.1290   &\textbf{18.2981} \\
    201807--201907 &   R-LSE        &0.3250    &0.0878   &\textbf{11.8815}  \\
                   &   LSE          &0.0117    &0.0001    &0.0168  \\
    201707--201807 &   R-LSE        &-0.2005    &0.0380    &\textbf{4.8738} \\
                   &  LSE        &-0.0523    &0.0027    &0.3385  \\
    201607--201707 &     R-LSE        &-0.7813    &0.3269   &\textbf{59.9758} \\
                   & LSE          &-0.1764    &0.0311    &3.9661  \\
    201507--201607 &    R-LSE        &-0.5840    &0.2115   &\textbf{33.1283}  \\
                   & LSE         & 0.0391    &0.0015    &0.1890  \\
   \bottomrule
  \hline
\end{tabular}
\end{table}

Table~\ref{table:SP500} shows that the model from the Robust-LSE
performs better than the one from the ordinary LSE fit according to
both  the index $R^2$ coefficient and the $F$-statistic of
goodness-of-fit.
 Furthermore beyond the 5 years reported in the table, we
  have also repeated the same comparison for all the past 20 years of
   the  S\&P500 Index: at 1\% level, $F$-statistic is 18 times significant
  for the model fitted with the Robust-LSE, while it is the case for
  one only model fitted with the ordinary LSE.

\section{Conclusion}
\label{sec:con}

In this study, a robust liner regression model under
both mean and variance uncertainty in the response variable  is
investigated.
We use a G-normal distribution to represent the variance uncertainty,
and another  nonlinear random variable for the mean
uncertainty. These nonlinear random variables in fact encompass an
infinite family of distributions for the response variable, instead of
a single distribution in the classical regression model.
For  a given estimation loss criterion,
two estimation strategies, namely the min-max and the min-min
strategies are introduced for estimating the regression parameter. 
The theory of sublinear expectation allows us to characterize the
optimal parameters for the two estimation strategies.
By considering the square loss function, the method leads to 
the robust (upper and lower) least squares estimators that
capture the maximum volatility and minimum volatility in the
response variable.
Under mild conditions on the data generation process, the consistency
of the estimators for both the regression parameter and the parameters
of mean and variance uncertainty is established.
These theoretical results are confirmed by simulation experiments. The
usefulness of the approach is assessed favorably in three applications
in comparison to the existing regression methods
including the ordinary LSE and a benchmark robust regression
estimator.

Further investigation of the proposed method would include more
extensive real data analysis. It is also worth researching on
alternative data generation process for the
general  distribution-uncertain regression model~\eqref{eq:line}.

\appendix

\section{Preliminaries from the sublinear expect ion theory}
\label{sec:prelim}

In the following, we introduce the sublinear expectation theory which
is used to describe the infinite family of distributions. We suppose
that there are an infinite family of probabilities
$\{P_{\theta}\}_{\theta\in \Theta}$ behind the error $\varepsilon$, and the related distribution is defined
as $F_{\theta}(z)=P_{\theta}(\varepsilon\leq z),\ z\in
\mathbb{R},\ \theta\in \Theta$, where $\Theta$ is a given set. Based
on the given infinite family of probabilities
$\{P_{\theta}\}_{\theta\in \Theta}$, we introduce the representation
results of a sublinear expectation $\mathbb{E}[\cdot]$, which is defined on a linear space $\mathcal{H}$ of
real valued functions on $\Omega$. A sublinear expectation
$\mathbb{E}[\cdot]:\ \mathcal{H}\to \mathbb{R}$ satisfies, for $X,Y\in \mathcal{H}$,

(i).\ \ \ $\mathbb{E}[X]\leq \mathbb{E}[Y],\ X\leq Y$;

(ii). \ $\mathbb{E}[c]=c,\ c\in \mathbb{R}$;

(iii). $\mathbb{E}[X+Y]\leq \mathbb{E}[X]+\mathbb{E}[Y]$;

(iv). $\mathbb{E}[\lambda X]=\lambda \mathbb{E}[X],\ \lambda\geq 0$.

The next result represents a sublinear expectation $\mathbb{E}[\cdot]$
as a supremum over a family of classical linear expectations.

\begin{theorem} \label{the:rep-sub}
  \citep[Theorem 1.2.1]{Peng2019}
Let $\mathbb{E}[\cdot]$ be a sublinear expectation on $\mathcal{H}$. There exists an infinite family of linear expectation $\{{E}_{\theta},\ \theta\in \Theta\}$ such that
\begin{equation}\label{eq:repre}
\mathbb{E}[X]=\max_{\theta\in \Theta}E_{\theta}[X],\ X\in \mathcal{H}.
\end{equation}
\end{theorem}

Define the space  $C_{l.Lip}(\mathbb{R})$ of  functions
 $\phi(\cdot)$ which are  locally Lipschitz: for some positive
 constants $C$ and $k$ depending on $\phi$, 
$$
\left| \phi(x)-\phi(y)\right|\leq C(1+\left|x\right|^k+\left|y\right|^k)\left|x-y\right|,\ x,y\in \mathbb{R}.
$$
We have the following  nonlinear  central limit theorem.

\begin{theorem}\label{the:clt}
  \citep[Theorem 2.4.4]{Peng2019}
Let $\{Z_i\}_{i=1}^{\infty}$ be a sequence of real-valued random variables on a sublinear expectation $(\Omega,\mathcal{H},\mathbb{E}[\cdot])$. Further, let $Z_{i+1}$ and $Z_i$  be identically distributed and $Z_{i+1}$ is independent from $\{Z_1,Z_2,\cdots,Z_i\}$ for $i\geq 1$. In addition, we assume that
$$
\mathbb{E}[Z_1]=\mathbb{E}[-Z_1]=0,
$$
and $\mathbb{E}[\left|Z_1\right|^{2+\delta}]<\infty$ for some $\delta>0$.
Then, the sequence
$$
\bigg{\{}\frac{Z_1+Z_2+\cdots+Z_n}{\sqrt{n}}\bigg{\}}_{n=1}^{\infty}
$$
converges to a G-normally distributed random variable $Z$ under sublinear
expectation $\mathbb{E}[\cdot]$: that is, for  $\phi(\cdot)\in C_{l.Lip}(\mathbb{R})$,
$$
\lim_{n\to \infty}\mathbb{E}\bigg[\phi(\frac{Z_1+Z_2+\cdots+Z_n}{\sqrt{n}})\bigg]=\mathbb{E}[\phi(Z)].
$$
\end{theorem}

The impact of this nonlinear central limit theorem on statistics is as
follows. 
In parallel to the role of the normal distribution that appears in the
limit of a classical central  limit theorem, the nonlinear G-normal random
variable $Z$ that appears in this theorem can serve as a natural model
for measurement errors in the nonlinear expectation framework, that
is, when variables are subject not to a single distribution but to
potentially infinite many and unknown distributions.  In this paper,  we apply this
idea to  the measurement error $\varepsilon$  in a linear regression
model as a way to catch up with its distribution uncertainty.

In the following, we develop more details on the  G-normal
distribution.

\subsection{The G-normal distribution for variance uncertainty}
\label{sec:prelim-Gn}

In the following, we explicitly construct a random variable $Z_1$ which follows
the  G-normal distribution given in Theorem \ref{the:clt}. Recall the
infinite family of probabilities $\{P_{\theta}\}_{\theta\in \Theta}$
introduced in Section \ref{sec:rlr}.
Let $Z_1$ satisfies
$$
\mathbb{E}[Z_1]=-\mathbb{E}[-Z_1]=x.
$$
Since
$\mathbb{E}[\cdot]=\max_{\theta\in\Theta}E_{\theta}[\cdot]$, this
relationship means that
$$
\max_{\theta\in\Theta}E_{\theta}[Z_1]=\min_{\theta\in\Theta}E_{\theta}[Z_1]=x,
$$
that is,
the maximum mean and the minimum mean of $Z_1$ over $\theta\in\Theta$
are the same.  In other words, $Z_1$ has no uncertainty on its
mean. The expectations of  $Z_1$ under
$\{P_{\theta}\}_{\theta\in \Theta}$ are given by
 \begin{equation}
   \label{eq:E_th}
   {E}_{\theta}[\phi(Z_1)]=\int_{R}\phi(z)dF_{\theta}(z),
 \end{equation}
where $\phi\in C_{l.Lip}(\mathbb{R})$ is some criterion (test) function.

In general, it is difficult to calculate the sublinear expectation
$\mathbb{E}[\phi(Z_1)]$.
We construct a G-normal distribution using a partial differential equation. This is because the partial differential equation tool can help us to find the optimal parameter $\theta_{\phi}$ such that $\mathbb{E}[\phi(Z_1)]={E}_{\theta_{\phi}}[\phi(Z_1)]$ and to calculate the expectation $E_{\theta}[\phi(Z_1)]$ under linear expectation $E_{\theta_{\phi}}[\cdot]$.
\begin{assumption} \label{ass:SDE}
 Let us assume $\{Z_t\}_{0\leq t\leq 1}$ satisfies the following stochastic differential equation,
   $$
   dZ_t=\theta_tdB_t,\quad Z_0=0,
   $$
   under $P_{\theta}$, $\theta\in \Theta=L^2(\Omega\times[0,1],[\underline{\sigma},\overline{\sigma}])$, where $\Theta$ is the set of all progressively measurable processes taking value on $[\underline{\sigma},\overline{\sigma}]$.
\end{assumption}

The stochastic process $\{Z_t\}_{0\leq t\leq 1}$ in Assumption
\ref{ass:SDE} admits a time-varying variance for the given probability
measure $P_{\theta}$. Therefore, there are infinite many distributions
behind this process.
We define the distribution of  
$Z_1$ as the G-normal distribution
$\EN_G(0,[\underline{\sigma}^2,\overline{\sigma}^2])$.\footnote{Based
on Assumption \ref{ass:SDE}, we use
$\EN_G(0,[\underline{\sigma}^2,\overline{\sigma}^2])$ to represent the
infinite family of distributions $\{F_{\theta}\}_{\theta\in\Theta}$ behind the random variable $Z_1$.}
Therefore, for a given criterion function $\phi(\cdot)\in C_{l.Lip}(\mathbb{R})$, we have 
$$
\mathbb{E}[\phi(Z_t)]=\max_{\theta\in\in\Theta}{E}_{\theta}[\phi(Z_t)]
=\max_{\theta\in\Theta}E_{\theta}[\phi(\int_0^t\theta_s\mathrm{d}B_s)].
$$
Proposition 2.2.10 of \cite{Peng2019} showed that $ u(t,x)=\mathbb{E}[\phi(Z_t+x)]$
is the unique viscosity solution of the following partial differential equation:
\begin{equation}
   \label{eq:pde-1}
   \partial_tu(t,x)-G(\partial_{xx}^2u(t,x))=0,\ t> 0,\ x\in \mathbb{R},
 \end{equation}
with the initial condition $u(0,x)=\phi(x),\ x\in \mathbb{R}$, where the function $G(\cdot)$ is defined as
\begin{equation}
  \label{eq:G-function}
  G(a) = \frac12  \left( \overline{\sigma}^2 a^+ -
    \underline{\sigma}^2a^-   \right),
  \quad a^+=\max(a,0), \ \text{and} \ a^-=\max(-a,0).
\end{equation}
It should be noted that $u(1,0)=\mathbb{E}[\phi(Z_1)]$.  Using the
process $\{Z_t\}_{0\leq t\leq 1}$, we can calculate the
characteristics of the G-normal random variable $Z_1$ for a given
criterion  function $\phi(\cdot)$ under
the infinite family of distributions $\{F_\theta\}_{\theta\in\Theta}$.

%
%

%

%
\section{Proofs}\label{sec:proofs}

\subsection{Proof of Lemma \ref{lem:minmax}}
\label{proof-lem-minmax}

In the first step, we prove that,
$$
\min_{\beta\in \mathbb{R}^q} \max_{\theta\in \Theta}E_{\theta}[\phi(Y-\beta^{\top} X-\mu_{\theta})]
=\max_{\theta\in \Theta} \min_{\beta\in \mathbb{R}^q} E_{\theta}[\phi(Y-\beta^{\top} X-\mu_{\theta})].
$$
For any given $\beta\in \mathbb{R}^q$,  since $X$ is a deterministic
vector variable, from (\ref{eq:line}), $\varepsilon=Y-\beta^{\top}
X-\eta$ satisfies a G-normal distribution
$\EN_G(a,[\underline{\sigma}^2,\overline{\sigma}^2])$, where $a$ is a
constant, which depends on $\beta$. Note that by assumption,
$\phi(\cdot)$ is  convex. Let 
$$
u(t,x)=\frac{1}{\sqrt{2\pi{\overline{\sigma}}^2t}}\int_{-\infty}^{\infty}\phi(y+x)
e^{-\frac{\left(y-a\right)^2}{2{\overline{\sigma}}^2t}}\mathrm{d}y.
$$
Because the  equation (\ref{eq:pde-1}) admits a unique classical
solution, we can verify that $u(t,x)$ is this solution, 
with  initial condition $\lim_{t\to 0}u(t,x)=\phi(x)$. Thus, we can take $\theta(s)=\overline{\sigma},\ 0\leq s\leq 1$ such that
$$
u(1,x)=E_{\overline{\sigma}}[\phi(Y-\beta^{\top} X-\mu_{\overline{\sigma}})]=\mathbb{E}[\phi(Y-\beta^{\top} X-\eta)].
$$
By Theorem \ref{the:rep-sub}, we have
$$
\mathbb{E}[\phi(Y-\beta^{\top} X-\eta)]=\max_{\theta\in \Theta}E_{\theta}[\phi(Y-\beta^{\top} X-\mu_{\theta})],
$$
and thus
$$
E_{\overline{\sigma}}[\phi(Y-\beta^{\top} X-\mu_{\overline{\sigma}})]=\max_{\theta\in \Theta}E_{\theta}[\phi(Y-\beta^{\top} X-\mu_{\theta})].
$$
It follows that
$$
\min_{\beta\in \mathbb{R}^q}E_{\overline{\sigma}}[\phi(Y-\beta^{\top} X-\mu_{\overline{\sigma}})]=\min_{\beta\in \mathbb{R}^q} \max_{\theta\in \Theta}E_{\theta}[\phi(Y-\beta^{\top} X-\mu_{\theta})].
$$
Obviously,
$$
\min_{\beta\in \mathbb{R}^q}E_{\overline{\sigma}}[\phi(Y-\beta^{\top} X)-\mu_{\overline{\sigma}}]\leq \max_{\theta\in \Theta} \min_{\beta\in \mathbb{R}^q} E_{\theta}[\phi(Y-\beta^{\top} X-\mu_{\theta})],
$$
which implies that
$$
\min_{\beta\in \mathbb{R}^q} \max_{\theta\in \Theta}E_{\theta}[\phi(Y-\beta^{\top} X-\mu_{\theta})]\leq \max_{\theta\in \Theta} \min_{\beta\in \mathbb{R}^q} E_{\theta}[\phi(Y-\beta^{\top} X-\mu_{\theta})].
$$

On the other hand, it is easy to verify that 
$$
\min_{\beta\in \mathbb{R}^q} \max_{\theta\in \Theta}E_{\theta}[\phi(Y-\beta^{\top} X-\mu_{\theta})]\geq \max_{\theta\in \Theta} \min_{\beta\in \mathbb{R}^q} E_{\theta}[\phi(Y-\beta^{\top} X-\mu_{\theta})].
$$
Thus, we have
\begin{equation}\label{eq:minmax1}
\min_{\beta\in \mathbb{R}^q} \max_{\theta\in \Theta}E_{\theta}[\phi(Y-\beta^{\top} X-\mu_{\theta})]= \max_{\theta\in \Theta} \min_{\beta\in \mathbb{R}^q} E_{\theta}[\phi(Y-\beta^{\top} X-\mu_{\theta})],
\end{equation}
and
\begin{equation}\label{eq:minmax2}
\min_{\beta\in \mathbb{R}^q}E_{\overline{\sigma}}[\phi(Y-\beta^{\top} X-\mu_{\overline{\sigma}})]=\min_{\beta\in \mathbb{R}^q} \max_{\theta\in \Theta}E_{\theta}[\phi(Y-\beta^{\top} X-\mu_{\theta})].
\end{equation}

Similarly, we can obtain the "min-min=min-min" exchange rule:
\begin{equation}\label{eq:minmin1}
\min_{\beta\in \mathbb{R}^q} \min_{\theta\in \Theta}E_{\theta}[\phi(Y-\beta^{\top} X-\mu_{\theta})]= \min_{\theta\in \Theta} \min_{\beta\in \mathbb{R}^q} E_{\theta}[\phi(Y-\beta^{\top} X-\mu_{\theta})],
\end{equation}
and
\begin{equation}\label{eq:minmin2}
\min_{\beta\in \mathbb{R}^q}E_{\underline{\sigma}}[\phi(Y-\beta^{\top} X-\mu_{\underline{\sigma}})]=\min_{\beta\in \mathbb{R}^q} \min_{\theta\in \Theta}E_{\theta}[\phi(Y-\beta^{\top} X-\mu_{\theta})].
\end{equation}
This completes the proof. $\qquad \qquad \Box$

\subsection{Proof of Theorem \ref{the:minmax}}
\label{proof-the-minmax}

Note that $\phi(\cdot)$ is convex. By the representation results (\ref{eq:minmax2}) and (\ref{eq:minmin2}) of Lemma \ref{lem:minmax}, we have
$$
E_{\overline{\sigma}}[\phi(Y-\beta^{\top} X-\mu_{\overline{\sigma}})]= \max_{\theta\in \Theta}E_{\theta}[\phi(Y-\beta^{\top} X-\mu_{\theta})],
$$
and
$$
E_{\underline{\sigma}}[\phi(Y-\beta^{\top} X-\mu_{\underline{\sigma}})]= \min_{\theta\in \Theta}E_{\theta}[\phi(Y-\beta^{\top} X-\mu_{\theta})].
$$
This implies that
$$
\overline{\beta}^{*}(\phi)=\arg\min_{\beta\in \mathbb{R}^q} E_{\overline{\sigma}}[\phi(Y-\beta^{\top} X-\mu_{\overline{\sigma}})]=\arg\min_{\beta\in \mathbb{R}^q} \max_{\theta\in \Theta}E_{\theta}[\phi(Y-\beta^{\top} X-\mu_{\theta})],
$$
and
$$
\underline{\beta}^{*}(\phi)=\arg\min_{\beta\in \mathbb{R}^q}E_{\underline{\sigma}}[\phi(Y-\beta^{\top} X-\mu_{\underline{\sigma}})]= \arg\min_{\beta\in \mathbb{R}^q}\min_{\theta\in \Theta}E_{\theta}[\phi(Y-\beta^{\top} X-\mu_{\theta})].
$$
This completes the proof. $\qquad \qquad \Box$

\subsection{Proof of Theorem\ref{the-consistency}}
\label{proof-the-consistency}

For notation simplicity, we set
$A_j=\{(x_i,y_i)\}_{i=1+n_0(j-1)}^{n_0j},\ 1\leq j\leq K$, with
$\eta_j\in [\underline{\mu},\overline{\mu}]$, and $\varepsilon_j\in\EN(0,\sigma_j^2)$, $\sigma_j^2\in   [\underline{\sigma}^2,\overline{\sigma}^2]$, the total number of
samples is $T=n_0K$. For each group $A_j$, when $n\leq n_0$, there
exists integer $k_j$ such that the samples
$\{(x_i,y_i)\}_{i=k_j}^{k_j+n-1}\subset A_j$. Thus, we can find a
block  $B_l \{(x_i,y_i)\}_{i=l}^{l+n-1}$ belongs to the group of
$\{A_j\}_{j=1}^K$ with the smallest variance $\min_{1\leq j\leq K}\sigma_j^2$.

\noindent{(i)}.\quad
Recall that within the $l$th block with
data $B_l=\{(x_i,y_i)\}_{i=l}^{l+n-1}$,
using ordinary LSE as defined in Step 1-(i) of the
procedure, we obtain the ordinary LSE for the regression parameter
and block mean, namely $(\hat{\beta}_l,\hat{\mu}_l)$.
The mean squared error $\hat{\sigma}^2_l$  in the block is also easily
obtained.
Recall the observation (ii) given below \eqref{eq:u-variance}: if one
$B_l$ overlaps with two $A_j$ groups, say $A_{j_l}$ and $A_{j_{l+1}}$,
the mean squared error  $\hat{\sigma}^2_l$ will be larger than if
$B_l$ is contained in a single $A_j$ group.
Therefore, the minimum of these mean squared errors will be achieved
by one block $B_l$ which is included in a single $A_j$.
Thus, by Theorem \ref{the:u-lse} and Corollary \ref{coro:1}, when $n\leq n_0$,
$\hat{\beta}_{\hat{k}}$ and $\hat\sigma^2_{\hat k}=\min_{1\leq l\leq m}\hat{\sigma}^2_l$ are  consistent estimators for
$\beta$ and $\min_{1\leq j\leq K}\sigma^2_j$.
As the latter is assumed to converge to $\underline{\sigma}^2$ as $K\to \infty$,
we have $\hat\sigma^2_{\hat k}  \to \underline{\sigma}^2$ with
probability $1$ as $K\wedge n\to\infty$. In a similar manner, we
obtain  $\hat{\beta}_{\hat{k}} \to \beta$ with probability 1 as $K\wedge n\to\infty$.
\medskip

\noindent{(ii)}.\quad
In  (i) above, we have obtained the consistency of
$(\hat{\beta}_{\hat k},\hat{\sigma}^2_{\hat k})$ for the parameters $(\beta,\min_{1\leq l\leq K}\sigma_l^2)$.
Recall the observation (i) given below \eqref{eq:u-variance},
the estimators
$\hat{\underline{\mu}}=\min_{1\leq l\leq m} \tilde{\mu}_l$  and $\hat{\overline{\mu}}=\max_{1\leq l\leq m} \tilde{\mu}_l$
for minimum mean and maximum mean from the groups $(A_j)_{1\le j\le  K}$
are consistent, that is,
$\hat{\underline{\mu}}=\min_{1\leq l\leq m} \tilde{\mu}_l$  and $\hat{\overline{\mu}}=\max_{1\leq l\leq m} \tilde{\mu}_l$
converge almost surely to $\min_{1\leq j\leq K}\eta_j$ and $\max_{1\leq j\leq K}\eta_j$.
As by assumption, the latter values converge to
$(\underline{\mu},\overline{\mu})$ as $ K \to \infty$,
the strong consistency of $(\hat{\underline{\mu}},\hat{\overline{\mu}})$ for
$(\underline{\mu},\overline{\mu})$  is obtained.

\medskip


\noindent{(iii)}.\quad Similar to the arguments given in (i), the upper variance estimator
$\hat{\overline\sigma}^2$ given in \eqref{eq:u-variance} converge
to $\max_{1\leq j\leq m}{\sigma}^2_j$.
As by assumption, the latter is assumed to converge to $\overline{\sigma}^2$ as $K\to \infty$,
we have $ \hat{\overline\sigma}^2   \to \overline{\sigma}^2$ with probability $1$.

The proof is complete. $ \qquad \qquad \Box  $

\bibliography{gexp,robuststat}

\end{document}